\documentclass[11pt]{amsart}

\usepackage[USenglish]{babel}

\usepackage{amsmath,amsthm,amssymb,amscd,bbm,verbatim}
\usepackage{booktabs}
\usepackage{url}

\usepackage{MnSymbol}

\usepackage{tikz}
\usetikzlibrary{arrows,shapes.geometric,positioning,decorations.markings}

\usepackage[mathscr]{eucal}
\usepackage[normalem]{ulem}
\usepackage{latexsym,youngtab}
\usepackage{multirow}
\usepackage{epsfig}
\usepackage{parskip}
\usepackage[all,cmtip]{xy}

\usepackage[margin=2cm]{geometry}
\usepackage{color}
\newtheoremstyle{custom}% name
  {3pt}%      Space above
  {3pt}%      Space below
  {\slshape}%         Body font
  {}%         Indent amount (empty = no indent, \parindent = para indent)
  {\bfseries}% Thm head font
  {.}%        Punctuation after thm head
  { }%     Space after thm head: " " = normal interword space;
   {}%         Thm head spec (can be left empty, meaning `normal')
\theoremstyle{custom}
\newtheorem{theorem}{Theorem}[section]
\newtheorem{proposition}[theorem]{Proposition}

\newtheorem{proposition/definition}[theorem]{Proposition/Definition}
\newtheorem{lemma}[theorem]{Lemma}
\newtheorem{corollary}[theorem]{Corollary}

\theoremstyle{definition}
\newtheorem{definition}[theorem]{Definition}

\newtheorem{example}[theorem]{Example}

\newtheorem{question}[theorem]{Question}

\theoremstyle{remark}
\newtheorem{remark}[theorem]{Remark}

\newtheoremstyle{exercise}% name
  {3pt}%      Space above
  {6pt}%      Space below
  {}%         Body font
  {}%         Indent amount (empty = no indent, \parindent = para indent)
  {\bfseries}% Thm head font
  {:}%        Punctuation after thm head
  { }%     Space after thm head: " " = normal interword space;
   {}%         Thm head spec (can be left empty, meaning `normal')
\theoremstyle{exercise}
\newtheorem{exercise}[theorem]{Exercise}
\newtheoremstyle{exercises}% name
  {3pt}%      Space above
  {6pt}%      Space below
  {}%         Body font
  {}%         Indent amount (empty = no indent, \parindent = para indent)
  {\bfseries}% Thm head font
  {:}%        Punctuation after thm head
  {\newline}%     Space after thm head: " " = normal interword space;
   {}%         Thm head spec (can be left empty, meaning `normal')

\theoremstyle{exercise}
\newtheorem{exercises}[theorem]{Exercises}

\input epsf
\def\boxit#1{\vbox{\hrule height1pt\hbox{\vrule width1pt\kern3pt
  \vbox{\kern3pt#1\kern3pt}\kern3pt\vrule width1pt}\hrule height1pt}}

\def\trank{\text{rank}}

\def\BC{\mathbb C}\def\BS{\mathbb S}
\def\BH{\mathbb H}
\def\BP{\mathbb P}
\def\pp#1{\mathbb P^{#1}}

\def\fgl{\mathfrak g\mathfrak l}

\def\pp#1{{\mathbb P}^{#1}}
\def\tdim{{\rm dim}}

\def\ww{\wedge}

\def\inv{{}^{-1}}

\def\cB{{\mathcal B}}\def\cA{{\mathcal A}}

\def\cE{{\mathcal E}}
\def\cF{{\mathcal F}}

\def\cO{{\mathcal O}}

\def\11{\mathbf 1}

\def\fsl{{\mathfrak {sl}}}

\def\fso{{\mathfrak {so}}}

\def\fg{{\mathfrak g}}

\def\a{\alpha}

\def\b{\beta}

\def\s{\sigma}

\def\ot{{\mathord{ \otimes } }}
\def\op{{\mathord{\,\oplus }\,}}

\def\lra{{\mathord{\;\longrightarrow\;}}}
\def\ra{{\mathord{\;\rightarrow\;}}}

\def\La#1{\Lambda^{#1}}

\def\cV{{\underline{V}}}

\def\frak{\mathfrak}
\def\fgl{\frak g\frak l}\def\fsl{\frak s\frak l}

\def\op{\oplus}
\def\BH{\Bbb H}\def\BZ{\Bbb Z}

\def\tann{\text{Ann}\,}

\def\ep{\epsilon}
\def\op{\oplus}

\def\ul{\underline}

\def\s{\sigma}
\def\t{\tau}

\def\a{\alpha}
\def\b{\beta}

\def\fso{\frak{so}}

\def\ol{\overline}

\def\BP{\mathbb  P}
\def\BC{\mathbb  C}

\def\pp#1{\mathbb  P^{#1}}

\def\BS{\mathbb  S}

\def\ep{\epsilon}

\def\opc{\op\cdots\op}

\def\cQ{\mathcal  Q}

\def\fg{\mathfrak  g}

\def\hd{, \dotsc ,}

\def\inv{{}^{-1}}

\def\La#1{\Lambda^{#1}}

\def\pp#1{\mathbb  P^{#1}}

\def\ur{\underline{\mathbf{R}}}

\def\ra{\rightarrow}

\def\tdet{\operatorname{det}}

\def\ttrace{\operatorname{trace}}
\def\tend{\operatorname{End}}
\def\tim{\operatorname{Im}}
\def\tdim{\operatorname{dim}}
\def\tker{\operatorname{ker}}\def\lker{\operatorname{Lker}}\def\rker{\operatorname{Rker}}

\def\tlim{\lim}
\def\tmod{\operatorname{mod}}

\def\tmin{\operatorname{min}}
\def\tmax{\operatorname{max}}
\def\thom{\operatorname{Hom}}
\def\trank{\operatorname{rank}}

\def\ww{\wedge}

\def\bbb{{\mathbf{b}}}

\def\be{\begin{equation}}
\def\ene{\end{equation}}
\def\aaa{{\mathbf{a}}}
\def\bbb{{\mathbf{b}}}
\def\ccc{{\mathbf{c}}}

\def\trank{\mathbf{R}}

\newcommand{\isom}{\cong}

\def\rank{\operatorname{rank}}

\newcommand{\Id}{\operatorname{Id}}

\def\tzeros{{\rm Zeros}}

\def\Mn{M_{\langle \nnn \rangle}}\def\Mone{M_{\langle 1\rangle}}\def\Mtwo{M_{\langle 2 \rangle}}

\def\Mn{M_{\langle \nnn \rangle}}\def\Mone{M_{\langle 1\rangle}}

\def\Mtwo{M_{\langle 2\rangle}}

\def\cK{{\mathcal K}}

\def\trank{{\mathrm {rank}}}

\def\aaa{\mathbf{a}}
\def\bbb{\mathbf{b}}
\def\ccc{\mathbf{c}}

\def\nnn{\mathbf{n}}

\def\rig#1{\smash{ \mathop{\longrightarrow}
    \limits^{#1}}}

\def\bt{\bold t}

\def\Mn{M_{\langle \nnn \rangle}}\def\Mone{M_{\langle
1\rangle}}

\def\Mtwo{M_{\langle 2\rangle}}

\def\cK{{\mathcal K}}
\def\lam{\lambda}

\def\trank{{\mathrm {rank}}}

\def\aaa{{\bold a}} \def\ccc{{\bold c}}

\def\nnn{\bold n}

\def\rig#1{\smash{ \mathop{\longrightarrow}
    \limits^{#1}}}

\def\bx{{\bold x}}

\def\om{\omega}

\renewcommand{\a}{\alpha}
\renewcommand{\b}{\beta}

\renewcommand{\BC}{\mathbb{C}}

\renewcommand{\hat}[1]{\widehat{#1}}

 \renewcommand{\tilde}{\widetilde}

\def\bz{\mathbf z}
\def\cV{\mathcal V}
\def\tcoker{\mathrm{coker}}
 
\begin{document}

\author{Hang Huang and J. M. Landsberg}
\address{Department of Mathematics, Texas A\&M University, College Station, TX 77843-3368, USA}
 
\email[J.M. Landsberg]{jml@math.tamu.edu}
  \title{On Linear spaces of of matrices bounded rank}

\thanks{Huang   supported by NSF grant  DMS2302375, Landsberg   supported by NSF grant  AF-2203618}

\begin{abstract} Motivated by questions in theoretical computer science and quantum information theory, we study the classical
problem of determining linear spaces of matrices of bounded rank. Spaces of bounded rank three
were classified in 1983, and it has been a longstanding problem to classify spaces of bounded rank four.
Before our study, no non-classical example of such a space was known. We exhibit two non-classical examples
of such   spaces   and   give the full classification of basic spaces
of bounded rank four. There are exactly four such up to isomorphism. We also take steps to bring together the methods
of the linear algebra community and the algebraic geometry community used to study spaces of bounded rank.
\end{abstract}

\maketitle

\section{Introduction}
A linear subspace   $E\subset \BC^\bbb\ot \BC^\ccc$ is of {\it bounded rank} if for all
$e\in E$, $\trank(e)< \tmin\{\bbb,\ccc\}$. We say a space of bounded rank has {\it bounded
rank $r$} if $r$ is the maximal rank of an element of $E$. Fix $\aaa=\tdim E$.

\begin{example} \label{comprex}
Let $E$ be of the form $\begin{pmatrix} *& *\\ *& 0\end{pmatrix}$ where the blocking is $(k_1,\bbb-k_1)\times
(k_2,\ccc-k_2)$. Then $E$ has bounded rank (at most) $k_1+k_2$. These are called {\it compression spaces},
or more precisely $(k_1,k_2)$-compression spaces.
\end{example}

\begin{example}\label{skewex} Let $E=\La 2\BC^{2p+1}\subset \BC^{2p+1}\ot \BC^{2p+1}$ be the space
of $2p+1\times 2p+1$ skew symmetric matrices. Then $E$ is of bounded rank $2p$ as the rank
of a skew-symmetric matrix is always even.
\end{example}

\begin{example}\label{kosex} Let $E\subset \thom(E,\La 2 E)$ be given by $e\mapsto \{ v\mapsto e\ww v\}$, which
has bounded rank $\aaa-1$, the kernel of the map associated to $e$ is the line through $e$.
\end{example}

We will refer to these examples as the {\it classical} spaces of bounded rank.

The study of  spaces of bounded rank dates back at least to
  Flanders in 1962 \cite{MR136618}, who solved a conjecture on the maximal dimension of such a space posed by Marcus.
Once the problem is stated, the case $\aaa=2$ follows immediately from the  
Kronecker-Weierstrass normal form for pencils of matrices.

Atkinson-Lloyd \cite{MR587090} introduced the notion of {\it primitive  spaces} and proposed the classification
of such by rank. As was known for a long time, the case $r=1$ consists of  the $(0,1)$ and $(1,0)$ compression spaces.  They classified the $r=2$ case, where the only primitive space is Example \ref{skewex} with $p=1$.  In \cite{MR695915} Atkinson carried out the classification for $r= 3$, the only primitive examples are Example \ref{kosex}  and its
projections.   In particular, there are no non-classical
examples of   spaces of bounded rank when $r\leq 3$. In the same paper he observed that if one allows $r$ to be large, there are 
\lq\lq many\rq\rq\  such spaces. This work is reviewed in \S\ref{ALsect}.

Sylvester \cite{MR840165} introduced language from geometry (vector bundles, Chern classes) to
study the more restrictive problem of spaces of constant rank $r$.

Eisenbud-Harris \cite{MR954659} independently introduced these tools, where for the general question
one must deal with sheaves rather than vector bundles. They proposed a refinement of the classification problem to  the study of
{\it basic spaces} (which in particular, disallows projections of primitive spaces) and stated that they were unaware of a non-classical example of a basic space
of bounded rank four. They also refined the observation that there are many such spaces for $r$ large
by observing that such spaces arise as matrices appearing in linear parts of the minimal free resolutions
of sufficiently general projective curves. In particular, there are nontrivial  moduli of basic spaces for $r$ large.
Their work is reviewed in \S\ref{EHsect}.

To our knowledge, the first example of a non-classical space of bounded rank is due to Westwick in 1996
\cite{MR1374258}, which is even of constant rank (namely constant rank $8$). Since then numerous explicit examples have been found.
See \cite{landsberg2022equivariant}  and \cite{ManRoi} for two recent contributions in the special case of constant rank.

\subsection*{Acknowledgements} We thank Austin Conner, Harm Derksen, David Eisenbud, Mark Green, Joe Harris, 
Laurent Manivel, Mihnea Popa,  
  Jerzy Weyman, and Derek Wu for useful conversations. Landsberg especially thanks Harris for inviting him to
Harvard where this project began to take shape. 

\section{Results}

Basic spaces are reviewed in \S\ref{basicsect}. All   spaces of bounded rank may be deduced from the basic spaces.

\begin{theorem}\label{mainthm}[Main Theorem] Up to isomorphism, there are exactly four
  basic spaces of bounded rank four:
\begin{enumerate}
\item[(I)]\label{ex1}   $E\isom \La 2 \BC^5\subset \BC^5\ot \BC^5$,

\item[(II)]\label{ex2} $E\isom  \BC^5\subset \thom(\BC^5, \La 2 \BC^5)$,

\item[(III)] \label{ourwinnerAat} 
$$
E=
\begin{pmatrix}
a_1& & & &-a_3&-a_5\\ & a_1& & & -a_4&-a_6\\ & & a_1& & a_2&0\\
& & & a_1& 0&a_2\\
 a_2& 0& a_3& a_5&0&0\\
0&a_2& a_4& a_6&0&0\end{pmatrix}\subset \BC^6\ot \BC^6,
$$

\item[(IV)]\label{ourwinnerBat}
$$
E=\begin{pmatrix}
 a_1 & a_2&   0 & 0 &-a_5&   -a_6    \\
   0  &  a_1 &    0  & 0& 0&  -a_5 \\
  0 &   0 &    a_1& a_2&  a_3& a_4  \\
 0 &   0  &  0   &  a_1& 0 &   a_3 \\
  a_3&  a_4&    a_5&    a_6&  0 &0 \\
  0  &   a_3& 0&   a_5  &   0 &   0  
 \end{pmatrix}\subset \BC^6\ot \BC^6.
$$
  
\end{enumerate}
\end{theorem}

The proof is given in \S\ref{mainpf}. It  utilizes classical (Atkinson-Lloyd) and algebreo-geometric techniques. 

\begin{question} Eisenbud-Harris \cite{MR954659} show that once $r$ is large,   that the basic spaces of bounded rank $r$ 
have moduli.  What is the smallest $r$ where moduli appear?
\end{question}

Our new examples  are part  of  a general construction:
\begin{proposition}\label{blowupprop}
Let $E$ be an $\aaa$-dimensional  space of bounded rank $r$ $\bbb\times \ccc$ matrices
and let $F\subset \tend(\BC^k)$ be an $\aaa$-dimensional space of
commuting matrices. Take   bases of $E$ and $F$
and for each matrix of the basis of $E$ replace each entry with the matrix of the corresponding basis element
of $F$ multiplied by the value of the entry  to obtain a $\bbb k\times \ccc k$ matrix. Call the new space $\tilde E$. Then $\tilde E$ is a space of bounded rank at most $kr$.
%Equality holds if a rank $k$ linear form replaces a linear form times a rank $r$ matrix.
\end{proposition}

The proof of Proposition \ref{blowupprop} is given in \S\ref{ANLsect}. 

The space $\tilde E$ of Proposition \ref{blowupprop}  is interesting only if the $k\times k$ matrices are not simultaneously diagonalizable,
as otherwise it is isomorphic to the direct sum of $k$ copies of $E$.

The two new examples in the Main Theorem are when $E=\La 3 \BC^3\subset \BC^3\ot \BC^3$ is  put in Atkinson
normal form (see \S\ref{ANLsect}) and the three basis elements are respectively
replaced by 
$$
\begin{pmatrix} a_1 & 0\\ 0&a_1\end{pmatrix}, \begin{pmatrix} a_2 & 0\\ 0&a_2\end{pmatrix},
\begin{pmatrix} a_3 & a_5\\ a_4&a_6\end{pmatrix}, \ \ \ {\rm and}\ \ \ 
\begin{pmatrix} a_1 & a_2\\ 0&a_1\end{pmatrix}, \begin{pmatrix} a_3 & a_4\\ 0&a_3\end{pmatrix},
\begin{pmatrix} a_5 & a_6\\ 0&a_5\end{pmatrix}.
$$

\begin{remark} A more geometric construction that is closely related to 
the construction of Proposition \ref{blowupprop} is presented in Proposition \ref{kronconstr} below. 
\end{remark}

\begin{remark} Yet another method for constructing larger spaces of bounded rank from
smaller ones is given in \cite[\S 5]{landsberg2022equivariant}.
\end{remark}

We   prove several technical results about invariants of spaces of bounded rank.
We introduce {\it Atkinson invariants} which arise naturally in the study of Atkinson normal form
for spaces of bounded rank,  and show that they upper bound the first  Chern classes of the sheaves
introduced in \cite{MR954659} (Proposition \ref{c1vat}). We correct a misconception in \cite{MR954659} about these sheaves
(Proposition \ref{tequal}, Example \ref{notdual}). In the motivation described in \S\ref{background}, tensors
of minimal border rank play a role. We show that such tensors cannot give rise to interesting spaces
of corank one (Corollary \ref{nominbir}). 
 We exhibit new families of corank two spaces generalizing \eqref{ourwinnerAat}, \eqref{ourwinnerBat} in \S\ref{generalsect},
as well as several variants.

\subsection*{Overview} 
In \S\ref{tensect} we explain the correspondence between spaces of matrices and tensors, give background
information on the geometry of tensors,   analyze the geometry of the tensors associated to the
examples in Theorem \ref{mainthm}, and give several generalizations. In \S\ref{reviewsect} we review the work of Atkinson-Llyod, Atkinson,
and Eisenbud-Harris. In \S\ref{unitesect} we take steps to unite the linear algebra and algebraic geometry perspectives.
In \S\ref{mainpf} we prove Theorem \ref{mainthm}. In \S\ref{generalsect} we provide a few
additional examples of spaces of bounded rank.

\section{Interpretations and generalizations via the associated tensors}\label{tensect}
\subsection{Background}\label{background}
Throughout this article,  $A,B,C$ are complex vector spaces of dimensions $\aaa,\bbb,\ccc$.
We let $\{ a_i\}$, $\{ b_j\}$, $\{c_k\}$ respectively be bases of $A,B,C$.
 There is a  1-1 correspondence between   tensors $T\in A\ot B\ot C$,  up to $GL(A)\times GL(B)\times GL(C)$ equivalence
  where the
 induced map $T_A: A^*\ra B\ot C$ is injective,
and $\aaa$-dimensional linear subspaces of $B\ot C$, i.e., points of the Grassmanian $G(\aaa,B\ot C)$,
up to $GL(B)\times GL(C)$ equivalence, given by $T\mapsto T(A^*)$.

In the study of tensors, the tensors that are the least understood are those that are $1$-{\it degenerate}: those where  $T(A^*)\subset B\ot C$, $T(B^*)\subset A\ot C$ and $T(C^*)\subset A\ot B$ are of bounded rank.

A tensor $T\in A\ot B\ot C$ has (tensor)  {\it rank one} if it is of the form $T=a\ot b\ot c$. The  {\it rank} 
of $T$ is the smallest $r$ such that $T$ may be written as a sum of $r$ rank one elements, and the
{\it border rank} of $T$, $\ur(T)$ is the smallest $r$ such that $T$ is a limit of tensors of rank $r$. In geometric
language, $\ur(T)$ is the smallest $r$ such that $[T]\in \s_r(Seg(\BP A\times \BP B\times \BP C))
\subset \BP (A\ot B\ot C)$, the $r$-the secant variety of the Segre variety. If $T_A: A^*\ra B\ot C$
is injective, then $\ur(T)\geq \aaa$.
The tensor $T$ is called {\it concise} when all three such maps are injective.
 When   $T$ is concise, if it has border rank $\tmax\{\aaa,\bbb,\ccc\}$, one says that  $T$ has {\it minimal border rank}.
 
\subsection{Motivations for this project}
In addition to being a classical problem of interest in its own right, two motivations for this project are as follows:

\subsubsection*{Strassen's laser method}
The problem to either unblock Strassen's laser method  for upper bounding the exponent
of matrix multiplication \cite{MR882307},  or prove it has exhausted its utility,  began in \cite{MR3388238}.  The problem
is to find new tensors to use in the method that can prove better upper bounds on the exponent than  the big Coppersmith-Winograd tensor, 
or prove that no such exist. Minimal border rank $1$-degenerate tensors are a class that are not yet known
to have barriers to improving the method \cite{blser_et_al:LIPIcs:2020:12686}.

\subsubsection*{Geometric rank and cost v. value in quantum information theory} While the rank and border rank
of tensors have been studied for a long time, recently  attention has been paid to the notions
of subrank, border subrank, and the closely related notions of slice rank \cite{Taoblog} and geometric rank
\cite{kopparty2020geometric}. Spaces of bounded rank give rise  to tensors with degenerate geometric
rank. In this paper we expand upon the (at the time surprising) result of \cite{MR4471039}, where
upper bounds on  geometric rank imply lower bounds on border rank.

\subsection{Brief discussion of tensors and their geometry}
Given a tensor $T\in A\ot B\ot C$, define its (extended) symmetry group
$$
\hat G_T=\{ g\in GL(A)\times GL(B)\times GL(C) \mid g\cdot T=T\}
$$
and the corresponding symmetry Lie algebra $\hat \fg_T$. 
The actual symmetry group $G_T$ has dimension two less as the
map $GL(A)\times GL(B)\times GL(C)\ra GL(A\ot B\ot C)$ has a two dimensional kernel. 

Several important (classes of)  tensors are:
\begin{itemize}
\item The unit tensor: $\Mone^{\op m}\in \BC^m\ot \BC^m\ot \BC^m=A\ot B\ot C$, where
$\Mone^{\op m}=\sum_{j=1}^m a_j\ot b_j\ot c_j$.
\item The $W$-state, also known as a general tangent vector to the Segre variety:
$W=a_1\ot b_1\ot c_2+a_1\ot b_2\ot c_1+a_2\ot b_1\ot c_1\in \BC^2\ot \BC^2\ot \BC^2$.
\item  Given an algebra $\cA$, one may form its structure tensor $T_{\cA}$ obtained from the bilinear
 map $\cA\times \cA\ra \cA$ given by multiplication.
\item The matrix multiplication tensor $\Mn\in \BC^{n^2}\ot \BC^{n^2}\ot \BC^{n^2}$ is the special
case of $T_{\cA}$ when $\cA$ is the algebra of $n\times n$ matrices.
\item Let $G\subset GL(A)\times GL(B)\times GL(C)$ be a subgroup and let
$L\subset A\ot B\ot C$ be a trivial $G$-submodule. Any
nonzero $T\in L$  is a $G$-invariant tensor, i.e., $G\subseteq \hat G_T$. Such are particularly interesting
when $G$ is large. However,  already when $G$ is a regular one-dimensional torus these tensors
(called {\it tight tensors} in this case)
can have interesting properties. Strassen conjectures
that such tensors have minimal asymptotic rank. See, e.g.,  \cite{landsberg2019finding,MR4204580,MR1341854}. 
\item A special case is when $A=B=V$, $C=\La 2V^*$, $G=SL(V)$ and $T\in \thom(\La 2V, V\ot V)$
is the inclusion map. When $\tdim V$ is odd, the corresponding tensor realizes two distinct
spaces of bounded rank, which are cases (I), (II) of the main theorem when $\tdim V=5$.
When $\tdim V$ is even,  the generalization of  \eqref{ex2} is still of bounded rank.
\item A  special case of the previous is when $\tdim V=3$,  as then $\La 2V^*\simeq V$ and one
obtains   $T_{skewcw,2}\in \La 3 \BC^3\subset \BC^3\ot \BC^3\ot \BC^3$.
\end{itemize}

Another motivation for this project was to find new classes of tensors of interest. In what follows
we show that the   tensors (III), (IV)  of Theorem \ref{mainthm} have many interesting properties and generalizations.

Given $T\in A\ot B\ot C$ and $T'\in A'\ot B'\ot C'$, their {\it Kronecker product} is just their tensor
product considered as a three-way tensor: $T\boxtimes T'\in
(A\ot A')\ot (B\ot B')\ot (C\ot C')$.

\subsection{A geometric interpretation of the tensor associated to Case (IV)  and a tensor variant of Proposition \ref{blowupprop}}
 Case (IV)  is a special case of the following construction:
 
 \begin{proposition}\label{kronconstr} Let $T\in A\ot B\ot C$ be such that $T(A^*)$ has bounded
 rank $r$ and let $T'\in A'\ot B'\ot C'$ with $\aaa'= \bbb'=\ccc'=m$ be of minimal border rank.
 Then $T\boxtimes T'(A^*\ot {A'}^*)$ has bounded rank at most $rm$.
 \end{proposition}
 
 \begin{proof} The proposition clearly holds when $T'=\Mone^{\op m}$ has rank $m$, as then one just obtains the
 sum of $m$ copies of $T$.   Being of bounded rank is a Zariski closed condition and all 
 concise minimal border rank $m$ tensors are degenerations of $\Mone^{\op m}$, i.e.,
 in $\ol{GL_m\times GL_m\times GL_m\cdot \Mone^{\op m}}$.
 \end{proof}
 
 Case (IV)  is the special case $T=T_{skewcw,2}$ and $T'=W$.
 
 Neither construction strictly contains the other: Case (III)
 does not arise from the construction of Proposition \ref{kronconstr} and a minimal border
 rank $1$-degenerate tensor does not in general correspond to a space of commuting matrices.

 \subsection{A geometric interpretation of Case (III)  and generalizations}
 The matrix multiplication tensor may be defined as follows: Let $U,V,W$ be vector spaces and let
 $A=U^*\ot V $, $B=V^*\ot W $, $C=W^*\ot U $.
 Then the (possibly rectangular) matrix multiplication tensor is 
 the (up to scale) unique $GL(U)\times GL(V)\times GL(W)$-tensor
 in $A\ot B\ot C$, namely
 $M_{\langle U,V,W\rangle}:=\Id_U\ot \Id_V\ot \Id_W$. In bases we may write $\Mtwo$ as
 $$
 \Mtwo(A^*)=\begin{pmatrix} x^1_1&x^1_2& & \\ x^2_1&x^2_2& &   \\ 
 & & x^1_1&x^1_2  \\ & & x^2_1&x^2_2 \end{pmatrix}.
 $$

 Consider the following construction which may be considered as an augmentation of
 the matrix multiplication tensor
 (see \S\ref{generalsect} for an even further generalization): let $U,V,W$ be even dimensional vector
 spaces equipped with symplectic forms $\om_U,\om_V,\om_W$. Let
 $A=U \ot V\op W$, $B=V \ot W\op U$, $C=W \ot U\op V$.
 Note that $A\ot B\ot C$ has a four dimensional space of $Sp(U)\times Sp(V)\times Sp(W)$ invariant tensors
 with basis $\om_U\ot \om_V\ot \om_W$ (which, identifying $U\simeq U^*$ etc. using
 the symplectic form,  is $\Id_U\ot \Id_V\ot \Id_W$), $\om_U\ot \om_V$, $\om_V\ot \om_W$, $\om_U\ot \om_W$.
 Here $Sp(U)=Sp(U,\om_U)$ is the symplectic group preserving  $\om_U$.
 Let 
 \be\label{ts}
 T_{UVW}=\om_U\ot \om_V\ot \om_W+ 
 \om_U\ot \om_W+\om_V\ot \om_W
 \ene
  so $G_{T_{UVW}}\supset Sp(U)\times Sp(V)\times Sp(W)$.
 Then when $\tdim U=\tdim V=\tdim W=2$ we obtain 
 Case (III). Explicitly, Case (III)  may be rewritten, setting
 $\{u_1,u_2\}$ a basis of $U\simeq U^*$  etc. and  $x^i_j=u_i\ot v_j$,
 $$
 \begin{pmatrix} x^1_1&x^1_2& & &-w_2&  \\ x^2_1&x^2_2& & &   &-w_2\\ 
 & & x^1_1&x^1_2&w_1& 
 \\ & & x^2_1&x^2_2 & &w_2\\
 w_1 & & w_2& & & \\
 & w_1 & & w_2& & 
 \end{pmatrix}.
 $$ 
 This construction shows in particular  that  Case (III)  has $\BZ_2$ symmetry 
 by exchanging the $B$ and $C$ factors. Note further that the space of matrices $T(B^*)$ is
 not of bounded rank. By the fundamental theorem of geometric rank \cite{kopparty2020geometric}
 $T(B^*)$ must have a highly nontransverse intersection with some $\s_r(Seg(\BP A\times \BP C))$.
 Indeed, it intersects $\s_2(Seg(\BP A\times \BP C))$ in a $\pp 1$.
 
 Recall that    the quaternion algebra is isomorphic to the algebra of $2\times 2$ matrices, so its structure
 tensor is $\Mtwo$.
 There is a six dimensional algebra $\BS$ that sits between the quaternions and the octonions, called
 the {\it sextonions} \cite{MR2204753}. We thank L. Manivel for observing that Case (III)
 is related to the sextonions.
 
 \begin{proposition} The tensor   $T_{UVW}$ of   \eqref{ts} when $\tdim U=\tdim V=\tdim W=2$  is the structure tensor of the sextonions, i.e.,  $T_{\BS}(A^*)$ is
 Case (III).
 \end{proposition}
 \begin{proof}
 This will be more transparent if we instead examine $T_{\BS}(C^*)$ (which is isomorphic to $T_{\BS}(B^*)$).
 Take bases $v_1\ot w_1,v_1\ot w_2,v_2\ot w_1,v_2\ot w_2, u_1,u_2$ of $B$
 and $u_1\ot v_2,u_2\ot v_2,u_1\ot v_1,u_2\ot v_1,w_1,w_2$ of $A$ and write $y_{ij}=u_i\ot w_j$.
 Then
 \be\label{T(B^*)}
 T(B^*)=
 \begin{pmatrix} y_{22} & -y_{21} & & & & \\
  -y_{12} & y_{11} & & & & \\
& &   -y_{22} & y_{21} & & \\
& &   y_{12} & -y_{11}  & & \\
& -v_2& &v_1& y_{22}&y_{12}\\
v_2& & -v_1& &y_{21}&y_{11}
\end{pmatrix}
\ene
On the other hand the multiplication in $\BS\simeq U^*\ot U\op  U$ is given by \cite[Def. 3.11]{MR2204753}
$(X,\mu)*(X',\mu')=(XX',(\ttrace(X)\Id -X)\mu'+X'\mu)$.
Writing the matrix representing the action of $(X,\mu)$ out in bases and permuting, we obtain the result.
\end{proof}

\subsection{Degeneracy loci and symmetry Lie algebras}
Let  $T(A^*)\subset B\ot C$ be of bounded rank $r$, where in this paragraph   we allow that possibly $r$ is full rank.
Let $\Sigma_A:=\{\a\in A^*\mid \trank(T(\a))<r\}$ denote  the degeneracy locus of the space.
 Let $G_{\Sigma_A}\subset GL(A)$ denote its symmetry group and $\fg_{\Sigma_A}\subset \fgl(A)$
 its symmetry Lie algebra.
 Let $\pi_A(\fg_T)\subset \fgl(A)$ denote the projection of $\fg_T$ onto the first factor.
 Then $\pi_A(\fg_T)\subset \fg_{\Sigma_A}$.
 Thus in general   one obtains that 
 $$
 \fg_T\subseteq \fg_{\Sigma_A}\op \fg_{\Sigma_B}\op \fg_{\Sigma_C}.
 $$

\medskip

  Case (II) is a space of constant
rank, and has been well-studied.
Case (I)  drops rank over the Grassmannian $G(2,5)\subset \BP (\La 2 \BC^5)$. (Here and below
we are interested in the projective geometry of the map so we work projectively - see \S\ref{EHsect} for more details.)
The symmetry group of $G(2,5)$ is $SL_5$ which is indeed the symmetry group of this tensor.

\medskip
 
 Case (IV)  is $\BZ_3$-invariant as it is the Kronecker product of two $\BZ_3$-invariant tensors,
namely $T_{skewcw,2}\boxtimes W$. Here $T_{skewcw,2}\in \La 3\BC^3$ and $W\in S^2\BC^2$.
One has $\fg_{T_{skewcw,2}}=\fsl_3$ and 
$$
\hat\fg_W=\left\{ \begin{pmatrix} \lam_1 &\lam_{12}\\ 0&-\mu_1-\nu_1\end{pmatrix}, \ 
 \begin{pmatrix} \mu_1 &\mu_{12}\\ 0&-\lam_1-\nu_1\end{pmatrix}, \ 
  \begin{pmatrix} \nu_1 &\nu_{12}\\ 0&-\mu_1-\lam_1\end{pmatrix}
  \mid \lam_{12}+\mu_{12}+\nu_{12}=0\right\}
$$
Here the degeneracy locus in each factor is just a linear space, e.g.
$\Sigma_A=\{\a_1=\a_3=\a_5=0\}$, so its Lie algebra in each factor is the parabolic
stabilizing the three plane, which has dimension $27$.
This is perhaps easier to see when we write the space as
$$
\begin{pmatrix} X& Y\\ Y& \end{pmatrix}
$$
where $X,Y$ are spaces of skew-symmetric $3\times 3$ matrices, and the degeneracy locus is when $Y=0$.
Explicitly we have the $21$-dimensional Lie algebra
$$
\hat \fg_{T_{skewcw,2}\boxtimes W}=\left\{ \begin{pmatrix} a&b\\ 0&a\end{pmatrix}, \ 
 \begin{pmatrix} a &b\\ 0&a\end{pmatrix}, \ 
  \begin{pmatrix} a &b\\ 0&a\end{pmatrix}
  \mid a,b \in \mathfrak{sl}_3 \right\} \oplus (\mathrm{Id} \boxtimes \hat\fg_W).
$$
The large increase over $8+5$ (where $5=\tdim \hat\fg_W$)
 is striking.  
 
 Super-additivity of dimension of symmetry Lie algebras under
 Kronecker product fails for generic tensors: the Kronecker product of two generic tensors
 will have no nontrivial symmetries. For unit tensors $\Mone^{\op k}\boxtimes \Mone^{\op \ell}=\Mone^{\op k\ell}$
 so the (extended symmetry) jump is to $2k\ell$ compared with an expected $2k+2\ell$.

 \begin{question} How to characterize situations where under Kronecker product, the dimension
 of the symmetry Lie algebra is super-additive?
 \end{question}

\medskip
  
 The symmetry Lie algebra of  Case (III) is $20$-dimensional, and it contains 
$\fsl_2^{\op 3}$. Explicitly, it is
$$
\hat\fg_{T_{\BS}}=\fgl(U)\times \fgl(V)\times \fgl(W)\op U\ot V\ot W
$$
where on $A=U\ot V\op W$, the action of $U\ot V\ot W$
is $u\ot v\ot w.(u'\ot v'+w')=\om_U(u,u')\om_V(v,v')w$, and on $B = V \ot W \op U$, the action is $u\ot v\ot w.(v'\ot w'+u')=\om_U(u,u')v\ot w$ and the action on $C$ is $u\ot v\ot w.(w'\ot u'+v')=\om_V(v,v')u\ot w$. 
To see $U\ot V\ot W\subset \hat \fg_{T_{\BS}}$, 
write out $T_{\BS}$ in bases
\begin{align*}
T_{\BS}&=(u_1\ot v_1)\ot (v_2\ot w_1)\ot (u_2\ot w_2)-(u_2\ot v_1)\ot (v_2\ot w_1)\ot (u_1\ot w_2)\\
&-(u_1\ot v_2)\ot (v_1\ot w_1)\ot (u_2\ot w_2)+(u_2\ot v_2)\ot (v_1\ot w_1)\ot (u_1\ot w_2)\\
&-(u_1\ot v_1)\ot (v_2\ot w_2)\ot (u_2\ot w_1)+(u_2\ot v_1)\ot (v_2\ot w_2)\ot (u_1\ot w_1)\\
&+(u_1\ot v_2)\ot (v_1\ot w_2)\ot (u_2\ot w_1)-(u_2\ot v_2)\ot (v_1\ot w_2)\ot (u_1\ot w_1)\\
&+w_1\ot u_1\ot (u_2\ot w_2)-w_2\ot u_1\ot (u_2\ot w_1)
-w_1\ot u_2\ot (u_1\ot w_2)+w_2\ot u_2\ot (u_1\ot w_1)\\
&+w_1\ot (v_1\ot w_2)\ot v_2-w_2\ot (v_1\ot w_1)\ot v_2
-w_1\ot (v_2\ot w_2)\ot v_1+w_2\ot (v_2\ot w_1)\ot v_1
\end{align*}
and apply  a basis vector $u_i\ot v_j\ot w_k$ to $T_{\BS}$.  One gets a sum of six terms that cancel in pairs.
Here the degeneracy locus in $\BP T_{\BS}(A^*)$ is the quadric surface $\{x^1_1x^2_2-x^1_2x^2_1\}\subset \pp 3\subset
\BP A^*$. The parabolic preserving this is $12+6$ dimensional  and by the above calculation, the nilpotent
part of the Lie algebra projected onto each factor   is isomorphic to $U\ot V\ot W$ which is contained in the parabolic.

In the general case where $U,V,W$ are $2k$-dimensional, one obtains a $(2k)^2+2k$ dimensional space
of bounded rank $(2k)^2$. That the space is of bounded rank is transparent from the Atkinson normal
form introduced below.

\begin{question} This construction gives rise to a series of algebras generalizing the sextonions.
What properties to these algebras have? Is there interesting representation theory associated to them?
(Using the sextonions, one may construct the Lie algebra $\frak e_{7\frac 12}$.)
\end{question}

\begin{remark} In general, the Lie algebra of derivations of an algebra $\cA$ sits inside
the two factor symmetry group of its structure tensor:  $\hat \fg_{T_{\cA},BC}:=\{ (Y,Z)\in \fgl(B)\times \fgl(C)\mid (Y,Z).T_{\cA}=0\}$.
The Lie algebra of derivations of the quaterions $\BH$ is $\fso(3)\isom \fsl_2$ and that of the 
sextonions is an $8$-dimensional algebra containing $\fsl_2$ as its semi-simple Levi factor.
\end{remark}

\subsection{Border ranks}
We thank Austin Conner for providing the decompositions in the results below:
\begin{proposition}\label{IIIbr}
 The tensor $T_{skewcw,2}\boxtimes W$ corresponding to Case (IV) has border rank nine.
The following is a $\BZ_3$-invariant border rank nine  decomposition:

\begin{align*}
T_{skewcw,2}\boxtimes W= &\\ 
\tlim_{t\ra 0}\frac 1{t^6} \BZ_3\cdot [
&(a_1-\frac 34 t^2a_3-\frac 34 t^3a_4)\ot
(-\frac 23 b_1+tb_2-t^3 b_4)\ot
(\frac 49 c_1-\frac 23 t c_2-\frac 23 t^2c_3+\frac 23 t^3 c_4+\frac 13 t^5 c_6)
\\
&+
t^2(a_1-\frac 12 ta_2-t^4a_5)\ot (-b_3+tb_4)
\ot (c_1+\frac 12 t c_2-\frac 12 t^2c_3+t^4c_5)\\
&+
(\frac 23 a_1-\frac 14 t^2 a_3-\frac 34 t^3 a_4+t^5 a_6)\ot
(\frac 23 b_1-tb_2+t^2b_3-\frac 12 t^5 b_6)\ot
(\frac 23 c_1+tc_2)].
\end{align*}
Here the $\BZ_3$ indicates cyclically permuting the factors to obtain nine terms.
 \end{proposition}
 \begin{proof} Verifying the decomposition is a routine calculation. The lower bound follows
 from computing the rank of the $p=1$ Koszul flattening (see, e.g., \cite[Ch. 2, \S 2.4]{MR3729273}).
 \end{proof}

 \begin{remark} For any two tensors $T,T'$  one has $\ur(T\boxtimes T')\leq \ur(T)\ur(T')$. 
 Proposition \ref{IIIbr} shows $9=\ur(T_{skewcw,2}\boxtimes W)<\ur(T_{skewcw,2})\ur( W)=10$.  Examples where strict submultiplicativity of border rank  under Kronecker product
 are of interest and also potentially useful for proving upper bounds on the exponent of matrix multiplication.
 In particular $T_{skewcw,2}$ could potentially be used to prove the exponent is two, and strict submultiplicativity
 of its Kronecker square has already been observed \cite{conner_harper_landsberg_2023}.
 We are currently examining which other tensors exhibit strict submultiplicativity under Kronecker product with $W$
 to potentially advance the laser method. 
 \end{remark}

 \begin{proposition} The tensor $T_{\BS}$ associated to Case (III)  has border rank $10$.
 \end{proposition}
 \begin{proof}  The upper bound comes from the following decomposition, which is easily verified:
 \begin{align*}
 T_{\BS}=  &\\ 
\tlim_{t\ra 0}\frac 1{t^5}[&
( a_1-t^{2}a_3+t^4a_6)\ot (tb_2+t^3b_3+t^4b_4-tb_5+ b_6)\ot (tc_4- c_6)\\
 +
&  (  a_1+t a_2)\ot (t^4b_4+ b_6)\ot (-t^3c_1+\frac 12 t^2c_3+c_6)\\
  +
 & (- a_1+t^{2}a_3)\ot (t^2b_1+tb_2+t^3b_3+b_6)\ot (-t^4c_2+tc_4-tc_5- c_6)
 \\
  +
 &  (-  a_1+t^4a_6)\ot ( tb_5-  b_6)\ot ( t^3c_1-t^4c_2-\frac 12 t^2c_3- c_6)]\\
+
\frac 1{t^4}[
&(- a_1+t a_2+t^3a_5)\ot b_5\ot (-\frac 12 tc_5+c_6)\\
+
& a_2\ot (tb_5+ b_6)\ot (t^3c_1+\frac 12t c_5- c_6)\\
+
 &(  a_1-t^{2}a_3+t^3a_5)\ot( tb_1+b_5)\ot (t^3c_1-t^2c_3+tc_4- c_6)
 \\
    -
& ( a_1-t^{2} a_3+t^3a_4)\ot ( tb_2+ b_6)\ot c_5]
 \\
 + \frac 1{t^3}[
&  ( a_1-t^3a_5)\ot (-tb_1-t^2b_3-\frac 12 b_5)\ot (2t^3c_2-tc_3+c_5)\\
+
&(- a_2+2t^2a_4)\ot (-t^3b_3+\frac 12  b_6)\ot (-2t^3c_2+tc_3+c_5)]. 
\end{align*}
To obtain the lower bound we use the border substitution method \cite{MR3633766,MR3842382}.
Note already from Strassen's commutation conditions applied to $T_{\BS}(B^*)$, after
changing bases so one basis element is the identity and using it to identify  $T_{\BS}(B^*)$
as a space of endomorphisms, the commutator has full rank which implies
$\ur(T_{\BS})\geq 9$. The border substitution method
says that if for every Borel fixed hyperplane the commutator still has full rank, one obtains
an improvement of one in the estimate. Here there are two Borel fixed hyperplanes, obtained
respectively by setting $v_2=0$ and $y_{22}=0$. In both cases the commutator is still of full rank.
 \end{proof}
 
\section{Review of previous work}\label{reviewsect}

\subsection{Work of Atkinson-Llyod and Atkinson \cite{MR587090,MR695915}}\label{ALsect}
Given $T\in A\ot B\ot C$ with $E=T(A^*)\subset B\ot C$ of bounded rank, to make
connections with existing literature we write, for $\a\in A^*$,  $\phi(\a)=\phi_T(\a): B^*\ra C$ for the
induced map. 

Consider   compression spaces  of bounded rank $r$ of Example \ref{comprex}:
they may be described invariantly as:  there exist $B'\subset B^*$ of codimension $k_1$,
$C'\subset C$ of dimension $k_2$, with  $ k_1+k_2=r$,  such that $T(A^*)(B')\subset C'$.
 This notion is symmetric as in this case $T(A^*)({C'}^\perp)\subset {B'}^\perp$. 
 
 Since compression spaces are \lq\lq understood\rq\rq, one would like to eliminate
 them from the study as well as spaces that are sums of a compression space with another
 space. Hence the following definition:
A bounded rank $r$  space is {\it imprimitive} if there exists $H\subset B^*$ such that $\phi|_H$ is
 of bounded rank $r-1$ or if there exists $H\subset C^*$ such that $\phi^\bt|_H$ is of bounded
 rank $r-1$. In this case the bounded rank condition is inherited from the smaller space.
 Otherwise it is {\it primitive}.

Thus  to classify spaces of bounded rank,  it suffices to classify the primitive ones.

Atkinson and Llyod \cite{MR587090}    showed that for primitive spaces of bounded rank $r$ in $B\ot C$, with $\bbb\leq \ccc$,
either $\bbb=r+1$ and $\ccc\leq \binom{r+1}2$ or 
  there exist $r_1,r_2\geq 2$ with $r_1+r_2=r$,
$\bbb\leq r_1+1+\binom {r_2+1}2$ and $\ccc\leq r_2+1+\binom{r_1+1}2$.

Thus to classify primitive  spaces of bounded rank $4$, it suffices to classify them in
$\BC^6\ot \BC^6$, and
$\BC^5\ot \BC^\ccc$,   with $5\leq \ccc\leq 10$. 

\subsection{Work  of Eisenbud and Harris   \cite{MR954659}}\label{EHsect} \label{basicsect}
Given $T\in A\ot B\ot C$, with $T(A^*)$ of bounded rank $r$, Eisenbud and Harris observed 
that one gets a map between vector bundles
\begin{align}\label{themap}
\ul B^* & \lra   \ul C\ot \cO_{\BP A^*}(1)\\
\nonumber \searrow & \ \   \ \ \swarrow\\
\nonumber & \BP A^* 
\end{align}
The notation is that $\ul C:=C\ot\cO_{\BP A^*}$ is the trivial vector bundle with fiber $C$.
Let $\cE=\cE_A$ denote the image sheaf and $\cF=\cF_A$ the image sheaf of the map with the roles of $B,C$ reversed.
Since $T(A^*)$ has bounded rank $r$, both sheaves   locally free off of a codimension at least two subset
because $\cE,\cF$ are torsion free. We slightly abuse notation, writing $\phi$ for the horizontal map \eqref{themap}.
 
They show that if $T(A^*)$ is primitive and $c_1(\cE)=1$, then
the space is obtained as a projection of the classical  Example \ref{kosex}.   I.e., $B\subseteq A$ and $C\subseteq \La 2 A^*$.

They assert $\cE^{**}\isom \cF^*(1)$ and by symmetry $\cF^{**}\isom \cE^*(1)$, but this is false
in general, see Example \ref{JPLex}.
Their assertion  would imply  $c_1(\cE_A)+c_1(\cF_A)= r $, which they
do not  assert but instead just assert $c_1(\cE_A)+c_1(\cF_A)\leq  r $, which we verify and identify the failure of
equality to hold (Proposition \ref{tequal}).

\begin{remark} When one has a space of constant rank the assertion is correct
as will be clear by Proposition \ref{tequal} and moreover in this case $\cE\isom \cF^*(1)$
as when $\cE$ is a vector bundle,  $\cE\isom \cE^{**}$.
\end{remark}

In particular, to classify the first open case of $r=4$, it suffices to understand when $c_1(\cE)=c_1(\cF)=2$.

They remark that the only basic $r=4$  case they know of with $c_1(\cE)=c_1(\cF)=2$
 is $B\isom C=\BC^5$ and $A=\La 2 B\subset B\ot B$.

We always have $B^*\subseteq H^0(\cE)$ and $C^*\subseteq H^0(\cF)$. Explicitly,  the first
inclusion is $\b\mapsto ([\a]\mapsto T(\a,\b)\ot \a^*$.

  Eisenbud and Harris point out that imprimitivity  is equivalent to $\cE$  or $\cF$  having a trivial summand.
 To see this, take a complement $L$ to $H\subset B^*$ and consider
 $\phi|_L: L\ra \cE$. It is injective as $\phi$ has bounded rank $r>r-1$, giving the desired splitting.
 
 Note that since $\ul{H^0(\cE)}$ surjects onto $\cE$, we can equally well have a trivial quotient
 of $\cE$, i.e., a surjection $\cE\ra \cO$ as this induces a surjection
 $H^0(\cE)\ra \cO$ which is just a linear map between vector spaces giving the splitting of $\cE$.
 But this in turn is equivalent to having an inclusion $\cO\ra \cE^*$, i.e., that
 $h^0(\cE^*)>0$. 
 
 In summary:
 
 \begin{proposition} $E$ is primitive if and only if   $h^0(\cE^*)=h^0(\cF^*)=0$.
 \end{proposition}

To further eliminate redundancies such as $E$ being a subspace of a primitive space or a projection
of a larger primitive space, Eisenbud and Harris defined a vector space of bounded rank $T(A^*)$ to be \emph{basic} if
\begin{enumerate}
    \item it is {\it strongly indecomposable}:  $\cE$ and $\cF$ are indecomposable. This is is a strengthening  of primitivity. 
    \item it is {\it unexpandable}:  $\varphi$ and $\varphi^*(1)$ are the linear parts of kernels of maps of graded free $\BC[A^*]$ modules. This says the space is not a projection of a space of bounded rank in some larger space
    of matrices.
    \item it is {\it unliftable}:  it is not a proper subspace of a family of the same rank in $\mathrm{Hom}(B^*,C)$.
\end{enumerate}
 
Draisma \cite{MR2268360} gives a  sufficient condition for a space to be unliftable 
(he calls unliftability  {\it rank criticality}): 
 for  $E\subset Hom(B^*,C)$  a linear space  of morphisms
of generic rank $r$, define  the space of {\it rank neutral directions} 
\begin{align*}RND(E):&=\{ X\in Hom(B^*,C), \; X(Ker(Y))\subset Im(Y) \; \forall Y\in E, rank(Y)=r\}\\
&=\bigcap_{Y\in E, rank(Y)=r} \hat T_Y\sigma_r(Seg(\BP B \times\BP  C)).
\end{align*}
Here $\sigma_r(Seg(\BP B \times\BP  C))$ denotes the variety of rank at most $r$ elements
and $\hat T_Y$ its affine tangent space at $Y$.

\begin{proposition}\label{draismaprop}  \cite{MR2268360}
  $RND(E)$ always contains $E$ and in case
of equality, $E$ is  rank critical.   
\end{proposition}

\section{Towards uniting the  linear algebra and algebraic geometry perspectives}\label{unitesect}

   \subsection{Atkinson normal form} \label{ANLsect}

\begin{lemma} \cite[Lemma 3]{MR695915}\label{atkinlem}
Let $E\subset \BC^\bbb\ot \BC^\ccc$ be of bounded rank $r$ and suppose we have
chosen bases so that some matrix $X_1\in E$ is of the form
$$
X_1=\begin{pmatrix} \Id_r & 0\\ 0 &0\end{pmatrix}  .
$$
Then for any $X\in E$ we have,  with the same blocking,
\be\label{ANF}
X=\begin{pmatrix} \bx&W\\ U&0\end{pmatrix}
\ene
and
$U\bx^kW=0$ for all $k\geq 0$.
\end{lemma}
Note that it is sufficient to check $U\bx^kW=0$ for $0\leq k\leq r-1$ because higher powers of $\bx$ may
be expressed as linear combinations of $\bx^0\hd \bx^{r-1}$.
\begin{proof}
    That $X$ has the bottom right block equal to zero follows by
    considering appropriate minors of $X + \varepsilon X_1$.
  Consider
$X_1+tX$ 
  and the size $r+1$ minor 
$$
0=\tdet\begin{pmatrix}\Id_r+ t\bx &tW_j\\ tU^k& 0
\end{pmatrix}
$$  
where $W_j$ is the $j$-th column of $W$ and $U^k$ is the $k$-th row of $U$.
Using the formula for a block determinant and assuming $t$ is sufficiently small,
\be\label{atseries}
\tdet\begin{pmatrix}\Id_r+ t\bx & tW_j\\ tU^k& 0
\end{pmatrix}=t^2 U^k(\Id_r+t\bx)\inv W_j
=t^2U^k(\sum_{s=0}^\infty (-1)^s t^s \bx^s)W_j.
\ene
The coefficient of each power of $t$ must vanish. Since this holds for all $j,k$,  we conclude.
\end{proof}

We add the observation that the normal form is sufficient as well:

\begin{proposition}\label{ANFconverse} Let $E\subset \BC^\bbb\ot \BC^\ccc$ admit a presentation in Atkinson normal
form \eqref{ANF}. Then $E$ is of bounded rank $r$.
\end{proposition}
\begin{proof} 
The normal form assures all size $r+1$ minors containing the upper-left $r\times r$ block are zero.
Say there were a size $r+1$ minor having rank $s$ in its $U$ block and $t$ in its $W$ block.
We may normalize 
$$U=\begin{pmatrix}\Id_s &0\\ 0&0\end{pmatrix}, W=\begin{pmatrix} 0&0\\ \Id_t &0\end{pmatrix},
$$
where the blockings are respectively $(s,r-s)\times (s,\bbb-r-s)$ and $(t,\ccc-r-t)\times (r-t,t)$. 
The $UW=0$ equation implies the first $s$ rows of $W$ and the last $t$ columns of $U$ are zero, so $s+t\leq r$.
Write $r+1=s+t+u$. There must be a size $u$ minor in the $\bx$ block contributing, and that
minor must be in the upper right $(r-s)\times (r-t)$ block of $\bx$. But the $U\bx W=0$ equation implies
that the $\bx$ block must be zero in the first $s$ rows of this upper right  block and the last $t$ columns of this block.
So the contribution must lie in a block of size $(r-s-t)\times (r-s-t)$. But $u=r-s-t+1$, a contradiction.
\end{proof}

\begin{proof}[Proof of Proposition \ref{blowupprop}]
Put $E$ in Atkinson normal form. Then $\tilde E$ will also be in Atkinson normal form when the linear form times
$\Id_r$ is replaced by a rank $k$ linear form. If not, just add in such a multiple of $\Id_r$ to put the space
in the normal form, as this will not effect the $\tilde U\tilde \bx^j \tilde W=0$ equations.
\end{proof}

Atkinson essentially observes that  the  equations $U\bx^jW=0$  may be encoded as
\be\label{amatrix}
\begin{pmatrix}
U \\ U \bx\\ U \bx^2\\ \vdots \\ U \bx^{r-1}
\end{pmatrix}
\begin{pmatrix} W &\bx W &\bx^2W &\cdots & \bx^{r-1}W 
\end{pmatrix}
=(0),
\ene
where $(0)$ is $(\ccc-r)\times (\bbb-r)$.
Call the first matrix on the left $At_L$ and call the second $At_R$.
Recall that if $X,Y$ are respectively $p\times r$ and $r\times q$ matrices, with $p,q\geq r$
and $\trank(XY)\leq s$ then $\tdim  \lker (X)+\tdim  \rker Y\geq r-s$, so
in our case we have $(r-\trank(X))+(r-\rank(Y))\geq r$,
i.e. $\trank(X)+\trank(Y)\leq r$, with $X$ the first matrix and $Y$ the second.

\begin{definition}\label{anumber} Notations as above. 
Assume bases have been chosen generically to have the normal form.
Write $at_L$ for the rank of $At_L$ and $at_R$ for the rank of $At_R$
  and call these the {\it Atkinson numbers} associated to $E$.
\end{definition}
In particular we have
\be\label{atinequal}
at_L+at_R\leq r.
\ene

 Say $E\subset B\ot C$ has bounded rank $r$ and we fix $X_1=\begin{pmatrix} \Id_r & 0\\ 0&0\end{pmatrix}$,
i.e., we take a general point of $E$. $X_1: B^*\ra C$ defines subspaces $\tim_{\thom(B^*,C)} (X_1)\subset C$ and $\tker_{\thom(B^*,C)}(X_1)^\perp=\tim_{\thom(C^*,B)}(X_1)\subset B$.
Call these subspaces $C_1, B_1$.
Let $X\in E$ be another general  element. Consider $X: B^*\ra C\tmod C_1$, and $X: C^*\ra B\tmod B_1$. 
The normal form implies that $X( B_1^\perp)\subseteq  C_1$ and $X( C_1^\perp)\subseteq  B_1$
The ranks
of these maps give two additional invariants.
In the coordinate expressions, these are respectively the maps given by $U$ and $W$.
We now allow $X$ to vary, and we ask what is the dimension of the subspaces
of $(B_1^\perp)^* \ot C_1$ and $(C_1^\perp)^*\ot  B_1$ that are filled out?
These give two additional invariants, call them $d_U,d_W$.
Let $d_{UW}$ be the dimension of the combined subspaces, in bases, the number
of variables appearing in $U,W$, combined.  We know by conciseness (to avoid trivialities) that 
$(B_1^\perp)^* \ot C_1$ and $(C_1^\perp)^*\ot  B_1$ respectively surject onto $C_1,B_1$
 which lower bounds their dimensions.
 
 \begin{proposition}\label{dUone} If $r=\ccc-1$ and  $d_U=1$ or $r=\bbb-1$ and $d_W=1$ then $E$ is imprimitive.
 More generally, for any $\bbb,\ccc,r$, if $U$ can be normalized to be zero except in one column or $W$ can be normalized
 to be zero except in one row, then $E$ is imprimitive.
 \end{proposition}
 \begin{proof}
 Say $r=\bbb-1$ and $d_W=1$. Then in Atkinson normal form we may normalize $W=(w_1,0\hd 0)^\bt$, so
 $$
 E=\begin{pmatrix} * & w_1\\ * &0\\ U&0\end{pmatrix}
 $$
 where the blocking is $(\bbb-1,1)\times (1,\bbb-2,\ccc-\bbb+1)$.
 If we strike out the first row the corank is unchanged.
 In the general case, the same argument holds when
 $$
 E=\begin{pmatrix} * & w_1&\cdots & w_k\\ * &0&\cdots &0\\ U&0&\cdots &0\end{pmatrix}.
 $$

 \end{proof}

 \begin{corollary}\label{nominbir} Assume $\aaa\leq\bbb= \ccc$. A concise tensor $T\in A\ot B\ot C$ such that $T(A^*)$ is
 corank one and primitive cannot have minimal border rank $\bbb$.
 \end{corollary}
 \begin{proof}
 By \cite[Prop. 3.3]{https://doi.org/10.48550/arxiv.2205.05713}, a tensor of minimal border rank
 as in the hypothesis must have $d_U=d_W=1$.
 \end{proof}
 
 We often slightly abuse notation and let $U,\bx, W$ denote the matrices of linear forms induced by
 $E$ rather than individual matrices.
 
Let $E$ be primitive and 
 let $\tann(W)$ denote the submodule of length $r$ row vectors with entries polynomials in the variables appearing
 in $W$ that annihilate $W$.
 
 \begin{proposition}\label{expandprop} Let $E\subset B\ot C$ be concise of bounded rank $r$ and primitive. If
 $\tann(W)$ is generated by row vectors of linear forms and  $at_L= \tdim(\tann(W)) = \bbb-r+1$, where $\tdim(\tann(W))$ is the dimension of the vector space spanned by minimal generators of $\tann(W)$,  then $E$ is expandable.
 \end{proposition}
 \begin{proof} Put $E$
 in Atkinson normal
  $\begin{pmatrix} \bx&W\\ U&0\end{pmatrix}$.
 Write $u_\a$ for the rows of $U$. Let 
 $$
 (At_L)_{prim}=\langle u_\a\bx^j\mid 1\leq\a\leq \bbb-r, j\geq 0\rangle \subseteq \tann(W).
 $$
 
Let $u\in \langle \textrm{minimal generators of } \tann(W)\rangle \backslash \langle u_\a \mid 1\leq\a\leq \bbb-r\rangle$.
Then
$\tilde E=\begin{pmatrix} \bx&W\\ U&0\\ u&0\end{pmatrix}$ is also
of bounded rank $r$ by Proposition \ref{ANFconverse} as it satisfies the conditions of Atkinson normal form.
\end{proof}

\subsection{Chern classes and Atkinson numbers}

\begin{proposition}\label{c1vat} Notations as in \S\ref{EHsect}, then $c_1(\cE) = at_R-c_1(\t_R)$ and $c_1(\cF) = at_L-c_1(\t_L)$,
where $\t_R,\t_L$ are torsion sheafs that are described explicitly in the proof. In particular,
$c_1(\cE)\leq at_R$ and $c_1(\cF)\leq at_L$.
\end{proposition}
\def\bss{B^*_1}\def\bsb{B^*_2}\def\css{C_1}\def\csb{C_2}
\def\ubss{\ul B^*_1}\def\ubsb{\ul B^*_2}\def\ucss{\ul C_1}\def\ucsb{\ul C_2}

\begin{proof}[Proof of Proposition \ref{c1vat}]
Put the matrix of  $\varphi:\underline{B}^* \to \underline{C}(1)$ in Atkinson normal form:
$$
\begin{pmatrix} \varphi_\bx & \varphi_W \\ \varphi_U & 0\end{pmatrix}.
$$
Over the point $[\a]\in \BP V^*$ where $\phi_{[\a]}=\begin{pmatrix} \Id_r&0\\0&0\end{pmatrix}$,
let $\bsb\subset B^*$ denote $\tker(\phi_{[\a]})$ and let $\bss$ be a complement in $B^*$.
In other words, $\ubss$ is 
  the trivial subbundle of $\ul B^*$ generated by the first $r$ basis vectors. 
Similarly write $C=\css\op \csb$. 
Note that we may identify
$\bss$ and $\css$.
 
We analyze $c_1(\cE)$. Consider the following maps:
\begin{align*}
\phi_i \quad \colon \quad \ubsb (-1) \oplus \ubsb (-2) \oplus \cdots \op \ubsb (-i)& \ra \ucss \\ 
 t_1\opc t_i&\mapsto  \varphi_W(t_1)+\varphi_\bx \varphi_W(t_2)+ \cdots +\varphi_\bx^{i-1}\varphi_W(t_i). 
 \end{align*}
 Note that $Im(\phi_i) = Im(\phi_{i+1})$ for all $i \geq r$. Denote the image of $\phi_r$ by $\cK_1$ and cokernel of $\phi_r$ by $\cQ_1$. Then $\rank \cK_1 = at_R$. Explicitly, we may  write $\cK_1$ as a subsheaf of $\ubss$ as follows:
 Let $\cV\subset \BP A^*$ be an affine open subset containing $[\a]$. Using the identification $B_1^*\isom C_1$,  
\[
\cK_1\mid_\cV = \{ s \in H^0( \bss\mid_\cV )  \  \mid \exists t_i \in  H^0(\ubsb(-i)\mid_\cV)  \text{\ such\ that\ }
 s = \varphi_W(t_1) + \varphi_\bx \varphi_W (t_2) + \ldots + \varphi_\bx^{r-1} \varphi_W (t_r)  \}.
\]

Since for all $k\geq 0$ and all  $t\in B_2^*$,  $\varphi_U\varphi_\bx^k\varphi_W(t)=0$,  
we obtain  $\varphi(\cK_1) \subseteq \cK_1(1)\subset \ucss(1)$  and $\rank \varphi\mid_{\cK_1} = \rank \cK_1 = at_R$.
We have the following commutative diagram that defines $\varphi^\cQ$:

\[
\xymatrix{
0 \ar[r] & \cK_1 \ar[r] \ar[d]^{\varphi\mid_{\cK_1}} & \ubss \oplus \ubsb \ar[r] \ar[d]^{\varphi} & \cQ_1 \oplus \ubsb \ar[r] \ar[d]^{\varphi^{\cQ}} & 0 \\
0 \ar[r] & \cK_1(1) \ar[r] & \ucss(1) \oplus \ucsb (1) \ar[r] & \cQ_1(1) \oplus \ucsb(1) \ar[r] & 0 
}
\]

Since $\varphi\mid_{\cK_1}$ is of full rank, $\tker(\varphi\mid_{\cK_1})$ is
torsion, but since it is a subsheaf of the  torsion free sheaf $\ul B_1^*$, it is zero,   and $\tcoker (\varphi \mid_{\cK_1})$ is a torsion sheaf. The first Chern classes in the short exact
sequence
$0\ra \cK_1\ra \cK_1(1)\ra \text{coker}(\varphi \mid_{\cK_1})\ra 0$ 
 show  
 $$c_1(\text{coker}(\varphi \mid_{\cK_1})) = \rank \cK_1 = at_R.
 $$
  Since coker$(\varphi \mid_{\cK_1})$ is torsion, any subsheaf of it has first Chern class no larger than $at_R$.

Claim: $c_1(\ker \varphi^{\cQ}) \geq 0$.

To see this, since $\varphi(\ubsb) = \varphi_W(\ubsb) \subseteq  \cK_1(1)$, we have $\ubsb \subseteq \ker \varphi^{\cQ}$. 
The kernel of the map $\varphi^{\cQ}\mid_{\cQ_1} \colon \cQ_1 \ra \cQ_1(1) \oplus \ucsb(1)$ is   torsion because 
$\varphi^{\cQ}\mid_{\cQ_1}$ is of full rank (since $\varphi\mid_{B_1^*}$ is of full rank). Hence $c_1(\ker \varphi^{\cQ}\mid_{\cQ_1}) \geq 0$. Together we have $c_1(\ker \varphi^{\cQ}) \geq 0$.

Hence we have $c_1(\tim( \varphi^{\cQ})) = c_1(\cQ_1\op \ubsb) - c_1(\ker \varphi^{\cQ})=c_1(\cQ_1) - c_1(\ker \varphi^{\cQ}) $.

We also have the following commutative diagram:

\xymatrix{
0 \ar[r] & \ker \varphi\mid_{\cK_1} \ar[r] \ar@{^{(}->}[d] & \ker\varphi \ar[r] \ar@{^{(}->}[d] & \ker \varphi^{\cQ} \ar[r] \ar@{^{(}->}[d] & \mathrm{coker} \varphi\mid_{\cK_1} \ar[r] \ar[d] & \mathrm{coker} \varphi \ar[r] \ar[d] & \mathrm{coker} \varphi^{\cQ} \ar[r] \ar[d] & 0 \\
0 \ar[r] & \cK_1 \ar[r] & \ubss \oplus \ubsb \ar[r] & \cQ_1 \oplus \ubsb \ar[r] & 0 \ar[r] & 0 \ar[r] & 0 \ar[r] & 0
}
whose rows are exact.

Define $\cA \colon = \ker(\mathrm{coker} \varphi\mid_{\cK_1} \to \mathrm{coker} \varphi)$ and let $\cB$ denote
the corresponding quotient. We have
\[
0 \longrightarrow \cA \longrightarrow \mathrm{coker} \varphi\mid_{\cK_1} \longrightarrow \cB \longrightarrow 0.
\]

Recall that  $\ker \varphi\mid_{\cK_1} = 0$. Consider the exact sequence
\[
0 \longrightarrow \ker \varphi \longrightarrow \ker \varphi^{\cQ} \longrightarrow \cA \longrightarrow 0.
\]
Thus  $c_1(\cE) = -c_1(\ker \varphi) = -c_1(\ker \varphi^{\cQ}) + c_1(\cA) = at_R - c_1(\tker(\varphi^\cQ)) - c_1(\cB)$.

The case of  $c_1(\cF)$ is similar.
\end{proof}

 \begin{example}\label{notdual} Let $T(A^*)\subset B\ot C$ correspond to a maximal dimension  $(k_1,k_2)$ 
 compression space of rank $r=k_1+k_2$ so that $\aaa=\bbb k_1+\ccc k_2-k_1k_2$.
 A standard way to write such a space is $\begin{pmatrix} *&*\\ * &0\end{pmatrix}$
 where the blocking is $(k_1,\bbb-k_1)\times (k_2,\ccc-k_2)$.
 Permute basis vectors to put it in Atkinson normal form
 \be\label{compranf}
 \begin{pmatrix} \bz_1&\bz_2&\bz_3\\ 0& \bz_4&0\\ 0 & \bz_5&0\end{pmatrix}
\ene
 where the blocking is $(k_1,k_2,\bbb-r)\times (k_1,k_2,\ccc-r)$.
 Then
 $$
 \begin{pmatrix}
U\\ U\bx\\ U\bx^2\\ \vdots \\ U\bx^{r-1}
\end{pmatrix}
=
\begin{pmatrix}
(0, \bz_5)\\
(0,\bz_5)\begin{pmatrix}\bz_1&\bz_2\\ 0 & \bz_3\end{pmatrix} \\
(0,\bz_5)\begin{pmatrix}\bz_1&\bz_2\\ 0 & \bz_3\end{pmatrix}^2 \\
\vdots
\end{pmatrix}=\begin{pmatrix}
(0, \bz_5)\\
(0,\bz_5\bz_3)  \\
\vdots
\end{pmatrix}
$$
Since the entries of $\bz_5$ are independent, we obtain $at_L=k_2$.
Since $\tdet(\bz_1)$ will factor out of the matrix of size $r$ minors from the first $r$ columns of \eqref{compranf},
but the remaining minors are independent because the entries of $\bz_4,\bz_5$
are independent, we also get $c_1(\cE)=k_2$.
Similarly $at_R=c_1(\cF)=k_1$.  If we specialize to a subspace of $\BP A^*$, then the
Atkinson numbers and Chern classes may drop.
 \end{example}

 \subsection{Chern class defects}
 Consider
 \be\label{QE}
0\ra \cE \ra \ul C(1)\ra \cQ_{\cE}\ra 0
\ene
and the analogous $\cQ_{\cF}$.  We sometimes write $\cQ=\cQ_\cE$ in what follows.
  
\begin{proposition} \label{tequal} Notations as in \S\ref{EHsect}, then 
$c_1(\cE)+c_1(\cF)=r-c_1(\mathcal{E}xt^1(\cQ_{\cE}, \mathcal{O}_{\mathbb{P} A^*}))$.
In particular $c_1(\mathcal{E}xt^1(\cQ_{\cE}, \mathcal{O}_{\mathbb{P} A^*}))=
c_1(\mathcal{E}xt^1(\cQ_{\cF}, \mathcal{O}_{\mathbb{P} A^*}))$. 
\end{proposition}

Note that  $c_1(\mathcal{E}xt^1(\cQ, \mathcal{O}_{\mathbb{P} A^*}))\geq 0$
because it is torsion.

We thank M. Popa for suggesting Proposition \ref{tequal} and   an outline of the proof.

\begin{proof}[Proof of Proposition \ref{tequal}]
 Dualize \eqref{QE} to get
\be\label{dualmap}
0\ra \cQ^*\ra \ul C^*(-1)\ra \cE^*\ra \mathcal{E}xt^1(\cQ, \mathcal{O}_{\mathbb{P} A^*}).
\ene
Now $\phi^*(1): \ul C^*\ra \ul B(1)$ is the map giving rise to $\cF$, so   $\phi^*: \ul C^*(-1)\ra \ul B$ has image  $\cF(-1)$.

Claim: we also have
  a map $\cE^*\ra \ul B$. 
To see this, take a resolution of $\cE$ by vector bundles
$$
\cV_2\ra \cV_1\ra \ul B^*\ra \cE\ra 0
$$
and dualize to get
$$
\cE^*\ra \ul B\ra \cV_1^*\ra \cV_2^*
$$
which  also gives
$$
0\ra \cQ^*\ra \ul C^*(-1)\rig{\phi^*}  \ul B\ra \cV_1^*\ra \cV_2^*
$$
which is {\it not} in general  exact at $\ul B$. By definition, the homology at that step is
$\cE xt^1(\cQ,\cO_{\BP A^*})$.

By \eqref{dualmap}, 
$\phi^*$ factors through the map $\cE^*\ra \ul B$.
We get an exact sequence
\be\label{FEmap}
0\ra \cF(-1)\ra \cE^*\ra \mathcal{E}xt^1(\cQ, \mathcal{O}_{\mathbb{P} A^*})\ra 0  
\ene
Thus  
$c_1(\cE^*)=c_1(\cF)-r+c_1(\mathcal{E}xt^1(\cQ, \mathcal{O}_{\mathbb{P} A^*}))$, i.e., 
$$
c_1(\cE)+c_1(\cF)=r-c_1(\mathcal{E}xt^1(\cQ, \mathcal{O}_{\mathbb{P} A^*})).
$$
 \end{proof}

 \begin{remark}
Dualizing \eqref{FEmap}, we see the  error in the assertion about $\cE^{**}$ and $\cF^*(1)$ is measured by 
\[
0 \ra \cE^{**} \ra \cF^*(1) \ra \mathcal{E}xt^1(\mathcal{E}xt^1(\cQ, \mathcal{O}_{\mathbb{P} A^*}),\mathcal{O}_{\mathbb{P} A^*})) \ra \mathcal{E}xt^1(\cE^*, \mathcal{O}_{\mathbb{P} A^*}) \ra \ldots.
\]
\end{remark}

\begin{example}\label{JPLex}

Consider the following  spaces from  \cite{https://doi.org/10.48550/arxiv.2205.05713}, which, as  tensors form a complete list of
concise,  $1$-degenerate tensors in $(\BC^5)^{\ot 3}$ of minimal border rank:

Represented as spaces of matrices, the tensors may be presented as:
\begin{align*}
 T_{\cO_{58}}&=
\begin{pmatrix} x_1& &x_2 &x_3 & x_5\\
x_5 & x_1&x_4 &-x_2 & \\
  & &x_1 & & \\
   & &-x_5 & x_1& \\
   & & &x_5  & \end{pmatrix}, 
   \ \
  T_{\cO_{57}} =
\begin{pmatrix} x_1& &x_2 &x_3 & x_5\\
 & x_1&x_4 &-x_2 & \\
  & &x_1 & & \\
   & & & x_1& \\
   & & &x_5  & \end{pmatrix}, 
\\
T_{\cO_{56}} &=
\begin{pmatrix} x_1& &x_2 &x_3 & x_5\\
    & x_1 +x_5 & &x_4 & \\
  & &x_1 & & \\
   & & & x_1& \\
   & & &x_5  & \end{pmatrix}, 
\ \ 
   T_{\cO_{55}}=
\begin{pmatrix} x_1& &x_2 &x_3 & x_5\\
 & x_1&  x_5  &x_4 & \\
 & &x_1 & & \\
   & & & x_1& \\
   & & &x_5  & \end{pmatrix}, \ \
   T_{\cO_{54}} =
\begin{pmatrix} x_1& &x_2 &x_3 & x_5\\
 & x_1& &x_4 & \\
  & &x_1 & & \\
   & & & x_1& \\
   & & &x_5  & \end{pmatrix}. 
 \end{align*} 

These are all bounded rank four.
The Chern classes  
$(c_1(\cE),c_1(\cF))$ equal the Atkinson numbers and
are respectively 
$(2,2),(1,1),(1,1),(1,1),(1,1)$. This is most easily computed via the Chern classes of the kernel bundles.

These are all compression spaces, as permuting the third and fifth columns makes
the lower $3\times 3$ block zero. They are degenerations of the general such compression space
over $\pp{15}$ where the kernel sheaves are isomorphic and move in a $Seg(\pp 1\times \pp 2)\subset \pp 5$,
which means $\cK=\cO_{\pp{15}}(-2)$. In particular the Chern classes can decrease under specialization.

Note  in the $(1,1)$ cases we have $c_1(\cE)=1$, $c_1(\cE^{**})=c_1(\cF^*(1))=4-1=3$, so
the Eisenbud-Harris assertions about $\cE$ and $\cF^*(1)$ fail in these cases.
 \end{example}

\section{Proof of the Main Theorem}\label{mainpf}

Outline of the proof:
 In \S\ref{zerosets}, we make general remarks on the $c_1(\cE)=2$ corank one case.
The $(5,5)$  case is treated  in \S\ref{55case}.
The $(5,m)$ cases, with $m>7$ are ruled out by  \cite[Thm. B]{MR695915}.
We refine the proof  of \cite[Thm. B]{MR695915} to rule out 
  the $(5,7)$ case in \S\ref{fiveseven}.  
We reformulate the Kronecker normal form for pencils of matrices and use
it to analyze linear annihilators of spaces of linear forms of $2\times r$ matrices in   \S\ref{kronformsect}. 
In \S\ref{66case} we use our analysis in \S\ref{kronformsect} to reduce to two new potentially
basic   spaces.
 The $(6,5)$ case is treated in \S\ref{65case}.
 In \S\ref{good Chern} we go to a slightly more general setting and observe the Chern classes
 in the two $(6,6)$ examples are indeed all $2$.
In   \S\ref{winnerisbasic} we prove the two $(6,6)$ examples are indeed basic.

\subsection{Zero sets and syzygies of spaces of quadrics}\label{zerosets}
Continuing the notation of \S\ref{EHsect}, write
$$
0\ra \cK_\cE\ra \ul B^*\ra C(1),
$$
and
$$
0\ra \cK_\cF\ra \ul C^*\ra B(1),
$$
which we dualize and twist to get
$$
\ul B^*\ra C(1)\ra \cK_\cF^*(1).
$$
Since the maps $B^*\ra C(1)$ agree we obtain
\be\label{sheafseqa}
0\ra \cK_\cE\ra \ul B^*\ra C(1) \ra \cK_\cF^*(1),
\ene
which is a complex but not in general exact.
Write the maps as $d_3,d_2,d_1$ respectively.

Recall the dictionary between sheaves on projective space and graded modules. 
Write $R=\BC[A]$ for the ring of polynomials on $A^*$  and adopt the usual convention that $R(k)$ is $R$ with the grading
shifted by $k$.  In terms of modules,
     assuming  $c_1(\cE)=c_1(\cF)=2$,   the corank one
  $(m,m)$ case  becomes:

   \begin{equation*}
\xymatrix@C+2pc{
 \mathbb{F}_{\bullet}  \quad \colon  &
  R(-2) \ar[r]^{\ker \varphi}
  &
  R^m \ar[r]^{\varphi} & R^m(1) \ar[r]^{\ker \varphi^\bt} & R(3).
}
\end{equation*} 
Here the entries of  the vector  $\tker\varphi$ are quadrics, which, by hypothesis have no common factor, and
similarly for the entries of 
$\tker\varphi^\bt$.

In this section we reduce the corank one case with $c_1(\cE)=2$. We already
saw that the kernel map is given by a row vector of quadrics which have no common factor.
We then specialize to the basic  rank four case and show the zero set of the quadrics always has codimension
at least three.

For $B\subset S^2A^*$, let 
$K_{21}(B)$ be the space of linear syzygies, i.e., the kernel of the multiplication map $B\ot A^*\ra S^3A^*$.
We will be interested in the corank one space of matrices given by the tensor
inducing $A^*\subset B\ot K_{21}(B)^*=:B\ot C$.

 Let $B\subset S^2A^*$, let   $\bbb=\tdim B\geq 3$ (the case of $\bbb=2$ is handled by Kronecker normal form).
 
Let $Q_1,Q_2\in B$ be general and consider their zero set. Then either
\begin{enumerate}
\item\label{ca}  The zero set is a degree four codimension two irreducible complete intersection. In this case
intersecting with a third general quadric in $B$ will provide a codimension three set.

\item\label{cb}   The zero set is of pure codimension two, consisting
of  the union of a cone over the twisted cubic $v_3(\pp 1)$ and a linear space.

\item\label{cd} The zero set is codimension two and the codimension two component consists of two  irreducible quadrics, each in a hyperplane.
In this case the space spanned by $Q_1,Q_2$ is $\langle Q, \ell m\rangle$, where $Q$ is irreducible and $\ell,m\in A^*$,
and the two irreducible quadrics are  $\tzeros(Q,\ell)$, $\tzeros(Q,m)$.
Then, by the genericity hypothesis, $B=\langle Q, \ell F\rangle$ where $\tdim F= \bbb-1$ and $\tdim K_{21}(B)= \binom{\bbb-1}2$.
\item\label{ce}  The zero set is codimension two and the codimension two components consist  of linear spaces.
Then $B=\langle \ell_1F_1,\ell_2F_2\rangle$, for some linear subspaces $F_1,F_2\subset  A^*$.

\item\label{cc}  The zero set has codimension one: then $B=\ell F$ where $\ell\in A^*$ and $F\subset A^*$ and then
$\tdim K_{21}(B)= \binom \bbb 2$.
\end{enumerate}

The only case requiring explanation is \eqref{cb}, which follows from the classification of varieties of
minimal degree, see, e.g., \cite{MR927946}, which in particular says that a degree three irreducible
variety of codimension two is a cone over the twisted cubic.

We now consider what possible basic spaces of bounded rank $\bbb-1$ in $B\ot C=B\ot K_{21}(B)$
could arise with kernel of dimension one given by the quadrics in $B$. We assume $\bbb>4$.
We claim only case \eqref{ca} can potentially occur:
Case \eqref{cb} cannot occur because in this case $\tdim K_{21}(B)=2$.
The remaining cases are ruled out because in each case there will be vectors of linear forms in 
the kernel.

 \subsection{Conclusion of the $(5,5)$ case}\label{55case}

\begin{lemma}\label{55out}
The only 
   basic space of $5 \times 5$ matrices of linear forms with $c_1(\cE) = c_1(\cF) = 2$ 
   is the  space of skew-symmetric matrices.
\end{lemma}

\begin{proof}
   Let $I$ denote the ideal generated by the quadrics in $\tker\varphi$.
By the discussion in \S\ref{zerosets},   $\mathrm{depth}\; I = \mathrm{codim}\; I  \geq 3$.

Set $I_k(\varphi)$ to be the ideal generated by the size $k$ minors of $\phi$. 
Using that  $\langle \La 4\phi\rangle=(\tker\phi)(\tker\phi^\bt)$, we have $\mathrm{gcd}(I_4(\varphi)) = 1$
as $c_1(\cE)=c_1(\cF)=2$. This implies $\mathrm{depth} \; I_4(\varphi) \geq 2$ and $\mathrm{depth} \; I_1(\ker \varphi) \geq 3$. This is enough for $\mathrm{coker} \; \varphi$ to be a first syzygy module, so that we can complete the resolution $\mathbb{F}_{\bullet}$ on the left:
    \begin{equation*}
\xymatrix@C+2pc{
 \mathbb{F}_{\bullet}  \quad \colon \quad 0 \ar[r] &
  R(-2) \ar[r]^{\ker \varphi}
  &
  R^5 \ar[r]^{\varphi} & R^5(1) \ar[r]^{\ker \varphi^\bt} & R(3)
}
\end{equation*}
By assumption, since $\varphi$ is primitive, the entries in $\ker \varphi^\bt$ are linearly independent. We are forced to have $\mathrm{depth} \; I_1(\ker \varphi^\bt) = \mathrm{codim} \; I_1(\ker \varphi^\bt) \geq 3$. By the
Buchsbaum-Eisenbud criterion for exactness \cite{MR314819},
 we conclude that both $\mathbb{F}_{\bullet}$ and its dual $\mathbb{F}_{\bullet}^*$ are exact. Hence $\mathbb{F}_{\bullet}$ is self dual (up to shift) and it defines a Gorenstein ideal of depth $3$. By the characterization of Gorenstein ideals of depth 3 \cite[Thm. 2.1(2)]{MR453723}, $\varphi$ is skew-symmetrizable.
\end{proof}

\subsection{$(r+1)\times m$ spaces of bounded rank $r$}\label{fiveseven}
A basis of the linear annihilator of  a vector of linear forms $(a_1\hd a_k)$ is
 $$
 \begin{pmatrix} -a_2\\a_1\\0\\ \vdots\\ \vdots\end{pmatrix},  \begin{pmatrix} -a_2\\ 0\\ a_1\\0\\ \vdots\end{pmatrix}, 
 \ldots , \begin{pmatrix} -a_k\\ 0\\ \vdots\\ 0 \\ a_1\end{pmatrix}, 
  \begin{pmatrix}0\\ -a_3\\a_2\\0\\ \vdots\\ \vdots\end{pmatrix},  \begin{pmatrix}0\\ -a_4\\ 0\\ a_2\\0\\ \vdots\end{pmatrix}, 
 \ldots , \begin{pmatrix}0\\ -a_k\\ 0\\ \vdots\\ 0 \\ a_2\end{pmatrix}, 
 \ldots ,
 \begin{pmatrix}0\\   \vdots\\ 0\\ -a_k \\ a_{k-1}\end{pmatrix}.
 $$
 In particular, it has dimension $\binom{k}2$.
 
 From this, one may recover   the Atkinson  result:
 \begin{theorem} \cite[Thm. B]{MR695915} Let $E\subset \BC^{r+1}\ot \BC^m$ be a primitive space of bounded rank
 $r$ and let $m>r +\binom{r-1}2$. Then $E$ is a specialization of $\BC^{r+1}\subset \thom(\BC^{r+1}, \La 2\BC^{r+1})$.
 \end{theorem}
 
 The idea of the proof is that assuming  $U$ has no entries equal to zero,
   there are enough elements of the linear annihilator of $U$ present so the $U\bx W=0$ equation forces
 $\bx$ to be a scalar times the identity. Writing out the matrix of linear forms of $V\ra \thom(V,\La 2 V)$ one
 sees this implies the space is a specialization of $V\ra \thom(V,\La 2 V)$.

 We now consider rank $4$ spaces in $\BC^5\ot \BC^7$:

 First assume   all entries of $U=(u_1,u_2,u_3,u_4)$ are linearly independent. In order to avoid
  $\bx$ being forced to be a linear form times the identity by the $U\bx W=0$ equation, $W$ is uniquely determined up to isomorphism,
 namely
 $$
 W=\begin{pmatrix}-u_2&-u_3&-u_4\\ u_1 & & \\ & u_1 & \\ & & u_1\end{pmatrix}.
 $$
 Here $u_1^2$ divides all size three minors, and one would obtain $c_1(\cF)=1$, so this case is eliminated.

 When   $U$ is three dimensional, normalize so $U=(0,u_2,u_3,u_4)$.
 Then if the first row of $W$ is zero, the full annihilator of $(u_2,u_3,u_4)$ must
 appear in $W$. Then using the $U\bx W=0$ equation we see that the first row of $\bx$ is zero
 except in the $(1,1)$ slot.  
 Otherwise the first row of $W$ contains a nonzero entry and  $U\bx W=0$
 gives the same conclusion for the first column of $\bx$. Striking the first row (resp. column) gives a space of bounded rank $r-1$ so the
 space is not primitive.    
 
 When $U$ is two-dimensional, write $U=(0,0,u_3,u_4)$, so we must have, after a change of basis
 $$
 W=\begin{pmatrix}  w^1_1&w^1_2&w^1_3\\ w^2_1&w^2_2&w^2_3 \\  & & -u_4\ep\\ &   &u_3 \ep \end{pmatrix}.
 $$
 with the upper left $2\times 2$ block of full rank to avoid a column of zeros after
 a change of basis and $\ep\in\{0,1\}$. Then   $U\bx W=0$ implies the lower left $2\times 2$
 block of $\bx$ is zero. Combining this with the first two columns of $W$ and below, we obtain
 a $4\times 3$ block of zeros and the space is compression.
 
 When $U$ is one-dimensional, the space is imprimitive by Proposition \ref{dUone}.

 We conclude:
 \begin{proposition}  Let $E\subset \BC^{5}\ot \BC^7$ be a primitive space of bounded rank
 $4$. Then $E$ is a specialization of $\BC^{5}\subset \thom(\BC^{5}, \La 2\BC^{5})$.
 \end{proposition}

\subsection{$2\times k$ spaces of linear forms and their linear annihilators}   \label{kronformsect}
Recall the  Kronecker normal form for pencils of matrices, i.e.,
tensors in $A\ot B\ot F$ with $\tdim F=2$: all are block matrices
with blocks
$$
L_k(F^*):=\begin{pmatrix} s&t& 0&\cdots &0\\
 0&s&t&0&\cdots &\\
 & & \ddots & &\\
 & & & s&t\end{pmatrix}, \ 
 L_k^\bt(F^*)=\begin{pmatrix} s&0& 0&\cdots &0\\
 t&s&0&0&\cdots &\\
 & & \ddots & &\\
 & & \ddots & &\\
 & & & t&s\end{pmatrix}, \
 Jor_{k,\lam}(F^*)= s\Id_k+tJ$$
 where the matrices are respectively $(k+1)\times k$, $k\times (k+1)$, $k\times k$ and in the last case $J$ is a single Jordan block
 with eigenvalue $\lam$.  
 Of particular note is
 $L_1(F^*)=(s,t)$ and $L_1^\bt(F^*)=\begin{pmatrix} s\\ t\end{pmatrix}$. We adopt the convention that
 the rows are indexed by $A$ and the columns by $B$.
 Rewritten as a linear subspace of $B\ot F$ these become:
 $$
 L_k(A^*)=\begin{pmatrix} a_1&a_2&\cdots & a_k&0\\ 0 & a_1& &\cdots & a_k\end{pmatrix}, \ 
L_k^\bt(A^*) = \begin{pmatrix} a_1&a_2&\cdots & a_{k }\\ a_2 &a_3  & \cdots & a_{k+1}\end{pmatrix},
$$
$$
 Jor_{k,\lam}(A^*)= \begin{pmatrix} a_1&a_2&\cdots & a_k \\ \lam a_1 & \lam a_2+a_1& \cdots & \lam a_k+a_{k-1}\end{pmatrix}.
 $$
 Note the special cases   
 $$L_1(A^*)=\begin{pmatrix} a&0\\0&a\end{pmatrix}, \ L_1^\bt(A^*)=\begin{pmatrix} a_1\\ a_2\end{pmatrix}, \ Jor_{1,0}(A^*)=\begin{pmatrix} a\\ 0\end{pmatrix}, \ Jor_{1,\lam}(A^*)=\begin{pmatrix} a\\ \lam a\end{pmatrix}, \
 Jor_{2,0}=\begin{pmatrix} a_1&a_2\\ 0&a_1\end{pmatrix}.
 $$
 
 Let  $0\leq u\leq v\leq f$,
 Write $A=A_1\opc A_{f}$. 
 The general form is,  
 $$
 (L_{k_1}(A_1^*)\hd L_{k_u}(A_u^*), L_{\ell_1}(A_{u+1}^*)\hd L_{\ell_v}(A_v^*),
 Jor_{i_1,\lam_{i_1}}(A_{v+1}^*)\hd Jor_{i_f,\lam_{i_f}}(A_f^*)).
 $$
 
 Now we study linear annihilators. First note that if a space consists of a single block, it has no linear annihilator.
 We next consider pairs of blocks. Label the linear forms in the first block with $a_j$'s and those
 in the second with $b_j$'s.
 
 An $L_k$ and an $L_q$, or
 an $L_k$ and a  $J_{q,0}$ with $q>1$, or
  a  $J_{k,0}$, $k>1$ and a  $J_{q,0}$, $q>1$: one gets a two dimensional space.  Respectively, 
   $$
 \begin{pmatrix} -b_1\\ 0\\ \vdots \\ 0 \\ a_1\\ 0\\ \vdots\\ 0\end{pmatrix},
  \begin{pmatrix}  0\\ \vdots \\ 0\\ -b_q  \\ 0\\ \vdots\\ 0\\ a_q\end{pmatrix} 
 ,\ \  {\rm and} \ \
 \begin{pmatrix} -b_1\\ 0\\ \vdots \\ 0 \\ a_1\\ 0\\ \vdots\\ 0\end{pmatrix},
  \begin{pmatrix}  -b_2\\ -b_1\\ 0  \\ \vdots \\ a_1  \\ a_2\\ 0 \\ \vdots\end{pmatrix}
 ,\ \  {\rm and} \ \ \begin{pmatrix} -b_1\\ 0\\ 0\\ \vdots \\ 0 \\ a_1\\ 0\\ 0\\ \vdots\\ 0\end{pmatrix},
  \begin{pmatrix}  -b_2\\ -b_1\\ 0\\ \vdots \\ 0 \\ a_2\\ a_1\\ 0 \\ \vdots\\ 0\end{pmatrix}.
    $$ 
    
 An $L_k$ with a $J_{1,0}$ has a one-dimensional linear annihilator. A $J_{k,0}$  (any $k$) with a $J_{1,0}$
 also has a one-dimensional linear annihilator.

An $L_1$ and an $L_q^\bt$ or a $Jor_{q,\lam}$: one gets a $q$ dimensional space. Respectively:
  $$
 \begin{pmatrix} -b_1\\ -b_2\\ 0\\ \vdots \\ 0 \\ a_1\\ 0\\ \vdots\\ 0\end{pmatrix} \hd 
  \begin{pmatrix} -b_q\\ -b_{q+1}\\ 0\\ \vdots \\ 0 \\ 0\\ 0\\ \vdots\\ a_1\end{pmatrix}
 ,\ \  {\rm and} \ \
 \begin{pmatrix} -b_1\\ -\lam b_1\\ 0\\ \vdots \\ 0 \\ a_1\\ 0\\ \vdots\\ 0\end{pmatrix},
   \begin{pmatrix} -b_2\\ -\lam b_2-b_{3}\\ 0\\ \vdots \\ 0 \\ 0\\ a_1 \\ 0\\ \vdots \end{pmatrix} \hd 
  \begin{pmatrix} -b_q\\ -\lam b_q-b_{q+1}\\ 0\\ \vdots \\ 0 \\ 0\\ 0\\ \vdots\\ a_1\end{pmatrix}.
$$

  All other pairs have no linear annihilator.
 
 \subsection{$(6,6)$-spaces}\label{66case}
 We now specialize to   spaces of $2\times 4$ matrices with linear forms as
 entries with at least a two-dimensional linear annihilator. The   discussion in \S\ref{kronformsect} implies the linear annihilator is 
 at most
 two-dimensional.  
 
 The same argument as in the $(5,7)$ case shows that $U$ cannot have any columns equal to zero.
 
 Going through partitions of four, we already know we need at least two parts.
  
 The case of $(3,1)$ is ruled out by our general analysis above. 
 
 We consider the $(2,2)$ and $(2,1,1)$ cases together.
 Thus  we have 
  $U=(u_1,u_2)$, where the variables appearing in $u_1$ are independent of those in $u_2$ and
 each has two columns, then $W=\begin{pmatrix} -u_2\\ u_1\end{pmatrix}$.
 By our analysis above, in order to have a two dimensional linear annihilator, 
 either $u_1=L_1=\mu\Id$ and $u_2$ can be anything,  or $u_1$ and $u_2$ are $J_{2,0}$'s.
 Note that in either case, $u_1,u_2$ must commute.
 We have the following possibilities:
 
 $u_1=L_1$ and $u_2$ anything, i.e., one of $L_1,J_{2,0},J_{2,\lam},L_1^\bt L_1^\bt, L_1^\bt J_{1,\lam},J_{1,\lam}J_{1,\lam'},
 J_{1,\lam}J_{1,0}
 L_1^\bt J_{1,0}, J_{1,0}J_{1,0}$ where $\lam,\lam'\neq 0$, or 
 
 $u_1=J_{2,0}$ and $u_2=J_{2,0}$.
 
 Note the only cases with $u_2$ non-invertible are $u_2=J_{1,0}J_{1,0}$ and $u_2=J_{1,\lam}J_{1,\lam}$ (same $\lam$ in both).

 Write $\bx=\begin{pmatrix} x_1& x_2\\ x_3&x_4\end{pmatrix}$ with each $x_j$ $2\times 2$, then
 the $U\bx W=0$ equation gives the equation  of $2\times 2$ matrices with cubic entries:
  $$
- u_1x_1u_2-u_2x_3u_2+u_1x_2u_1+ u_2x_4u_1=0
$$

Write $x_j=y_j+x_j(u_1)+x_j(u_2)$ where the variables in $y_j$ are independent of those in $u_1,u_2$
and $x_j(u_i)$ has entries linear in the entries of $u_i$. We immediately obtain
$y_2=x_2(u_1)=0$  and if $u_2\neq J_{1,0}J_{1,0}$, $y_3=x_3(u_2)=0$.

Note we are free to modify $x_1$ (resp. $x_3$) by a multiple of $u_1$ as long as we modify $x_2$
(resp. $x_4)$, by the same multiple of $u_2$, and $x_1$ (resp. $x_2$) by a multiple of $-u_2$
as long as we modify $x_3$ (resp. $x_4$) by the same multiple of $u_1$.

Consider the case $u_1=u_2=J_{2,0}$. Note that the centralizer (among matrices of linear forms) of a $J_{2,0}$ is spanned by
$t\Id$ and $J_{2,0}$. Separating the equations by the variables involved, we obtain
$y_1=y_4=t\Id + J_{2,0}$, and we may absorb the $t\Id$ into the $J_{2,0}$.
Taking into account the other homogeneities and our allowed modifications, we
obtain Case (IV).

Next consider the case $u_1=L_1$ and $u_2$ has   centralizer $t\Id$, which occurs when $u_2$
is any of  
$J_{2,\lam},L_1^\bt L_1^\bt, L_1^\bt J_{1,\lam},
 L_1^\bt J_{1,0}$. In all these cases we are reduced to $\bx= t\Id_4$, but the resulting spaces are all
 specializations of the case $u_2=L_1^\bt L_1^\bt$, and thus not basic unless $u_2=L_1^\bt L_1^\bt$, which is Case (III).
 
 Case $(u_1,u_2)=(L_1,L_1)$: the solutions to the homogeneous parts of the
 equations, after admissible modifications, give $y_1=y_4$, and all other terms in $\bx$ are zero.
 The resulting space is a permuted version of Case (III).

  Case $u_1=L_1$, $u_2=J_{2,0}$. This case gives a specialization of Case (IV).
 
 Case $u_1=L_1$, $u_2=J_{1,0}J_{1,0}$. The $U\bx W=0$ equation implies $x^2_2$ is the only nonzero
 entry in the second row, so this case is not primitive.

   The case $(1,1,1,1)$ is easily ruled out.

 \subsection{$(6,5)$ case}\label{65case}
 The same argument as in the $(5,7)$ case shows that $U$ cannot have any columns equal to zero. Thanks to Proposition \ref{expandprop}, the cases that might give basic spaces are when $U$ has a one-dimensional linear annihilator and $at_R=2$, so
 $ \bx W$ must not be a linear form times $W$, and thus $U$ has a primitive degree two annihilator.
 By the discussion in \S\ref{kronformsect}, there are four cases with a one-dimensional linear annihilator:
 $(L_2,J_{1,0})$, $(J_{3,0},J_{1,0})$,   $(J_{2,0},J_{1,0},L_1^\bt)$, and $(J_{1,0},J_{1,0},L_1^\bt,L_1^\bt)$ which we may respectively write as
 $$
 \begin{pmatrix}
 u_1&u_2&0&u_3\\ 0&u_1&u_2&0\end{pmatrix},\ 
  \begin{pmatrix}
 u_1&u_2&u_3&u_4\\ 0&u_1&u_2&0\end{pmatrix},\ 
 \begin{pmatrix}
 u_1&u_2&u_3&u_4\\ 0&u_1&0&u_5\end{pmatrix},\ 
  \begin{pmatrix}
 u_1&u_2&u_3&u_4\\ 0&0&u_5&u_6\end{pmatrix}.
 $$
 Note that the first is a degeneration of the second, and after permuting columns, the second is
 a degeneration of the third, and the third a degeneration of the fourth. To make the degeneration transparent, we rename the variables and write the spaces as follows
 $$
 \begin{pmatrix}
 u_1&u_2&u_3&0\\ 0&0&u_1&u_3\end{pmatrix},\ 
  \begin{pmatrix}
 u_1&u_2&u_3&u_4\\ 0&0&u_1&u_3\end{pmatrix},\ 
 \begin{pmatrix}
 u_1&u_2&u_3&u_4\\ 0&0&u_1&u_6\end{pmatrix},\ 
  \begin{pmatrix}
 u_1&u_2&u_3&u_4\\ 0&0&u_5&u_6\end{pmatrix}.
 $$
 In all the cases, $W = (-u_2,u_1,0,0)^\bt$.
 We computed by Macaulay2 that primitive degree two annihilator of $U$ is generated by $(0,-u_4 u_5 + u_3 u_6,-u_2 u_6, u_2 u_5)^\bt$ and $(-u_4 u_5 + u_3 u_6, 0 , -u_1 u_6, u_1 u_5)^\bt$, where the effect of degeneration does not change the space of primitive degree two annihilators. By analyzing the first entry of $\bx W$, we see that $\bx W$ has to be a linear form times $W$ since it cannot involve any primitive degree two annihilators, which is a contradiction.

 \subsection{Good Chern classes}\label{good Chern}
 By construction the Atkinson numbers in our examples are all $2$, here we see the Chern classes are as
 well. We work in a slightly more general setting where $u_1,u_2,x_1$ are  commuting $k\times k$ matrices,
 i.e., the general setting of Proposition \ref{blowupprop} applied to $\BC^3\subset\thom(\BC^3,\La 2\BC^3)$.
 Our space is
 $$
 \begin{pmatrix} x_1& & -u_2\\ & x_1&u_1\\ u_1&u_2& \end{pmatrix}
 $$
 Then the left kernel is the $k\times 3k$ matrix of linear forms $(-u_1,-u_2,x_1)$, and the
 right kernel is the $3k\times k$ matrix $\begin{pmatrix} u_2\\ -u_1\\ x_1\end{pmatrix}$.
 Since the blocks $u_1,u_2,x_1$ are each in different sets of variables, the matrix of size $k$ minors
 has no common factor and we conclude $c_1(\cE)=c_1(\cF)=k$.

\subsection{Proof that our two new examples are basic}\label{winnerisbasic}

\begin{proof} 
By the
Buchsbaum-Eisenbud criterion for exactness \cite{MR314819}, it is easy to see that the following complexes are exact as well as their duals:
\begin{equation*}
\xymatrix@C+2pc{
 0 \ar[r] &
  R^2 \ar[r]^{d_3}
  &
  R^6 \ar[r]^{Case (III)}
  & 
  R^6 \ar[r]^{d_1}
  &
  R^2.
}
\end{equation*}

\begin{equation*}
\xymatrix@C+2pc{
 0 \ar[r] &
  R^2 \ar[r]^{d_3^{'}}
  &
  R^6 \ar[r]^{Case (IV)}
  & 
  R^6 \ar[r]^{d_1^{'}}
  &
  R^2.
}
\end{equation*}
where $R = \mathbb{C}[a_1,\ldots,a_6]$, 
$$d_3 = \left( \begin{smallmatrix} a_3 & a_5 \\ a_4 & a_6 \\ -a_2 & 0 \\ 0 & a_2 \\ a_1 & 0 \\ 0 & a_1 \end{smallmatrix} \right),\    d_1 = \left( \begin{smallmatrix} -a_2 & 0 & -a_3 & -a_5 & a_1 & 0 \\ 0 & -a_2 & -a_4 & -a_6 & 0 & a_1
\end{smallmatrix}   \right), \ d_3^{'} = \left( \begin{smallmatrix} a_5 & a_6 \\ 0 & a_5 \\ -a_3 & -a_4 \\ 0 & -a_3 \\ a_1 & a_2 \\ 0 & a_1 \end{smallmatrix} \right), {\rm \  and \ }d_1^{'} = \left( \begin{smallmatrix} -a_3 & -a_4 & -a_5 & -a_6 & a_1 & a_2 \\ 0 & -a_3 & 0 & -a_5 & 0 & a_1
\end{smallmatrix}   \right).
$$
 This proves Case (III) and Case (IV)  are unexpandable.

The spaces Case (III) and Case (IV) are strongly indecomposable
if and only if  the  image sheaves of $d_1,d_1^{'}$ as well as the image sheaf of $d_3^\bt , (d_3^{'})^\bt$ are indecomposable. In the case corresponding to Case (III), we compute the  ideal generated by maximal minors of $d_3$ as well as $d_1$. They each have  $12$ minimal generators. However, if their image sheaf were decomposable, the number of minimal generators would be  at most $\binom{5}{2} = 10$. To see this, one obtains 
the most generators when the last column  of $d_3^\bt$ or $d_1$ is zero.  In the case corresponding to Case (IV), note that the entries of the first row vector of $d_1^{'}$ as well as the entries of the second column vector of $d_3^{'}$ are algebraically independent (even after any possible row operations). Hence up to row/column operations and $\mathrm{GL}(A)$ actions, we cannot have a zero entry in the first row of $d_1^{'}$ (respectively, the second column of $d_3^{'}$). Thus the image sheaves of $d_1^{'}$ and $d_3^{'}$ are indecomposable.

Finally, to show Case (III) and Case (IV) are unliftable, we used   code in Macaulay2 implementing  
Proposition \ref{draismaprop}. 
\end{proof}
 
\section{Additional Examples}\label{generalsect} We first give two corank two examples generalizing
the new rank four cases:
   If there are exactly two $L_1$ blocks and at most one  $Jor_{1,0}$ block, then  there is a  $2\bbb-6$ dimensional
 kernel, and the resulting space of bounded rank $\bbb-2$ size  $\bbb\times (2\bbb-6)$ matrices is a specialization of
 the case with $\bbb-4$ $L_1^\bt$'s, so it is sufficient consider that space.
 That is, set $p=\bbb-3$,the space is
 \be\label{twoL1}
 \begin{pmatrix} a_1&0&a_2&0&a_3&a_5&\cdots & a_{2p-1}\\
0&a_1&0&a_2&a_4&a_6& \cdots & a_{2p}\end{pmatrix}
\ene
and after permuting rows of the kernel matrix to make it in Atkinson normal form one obtains
 \be\label{twoL1bndd}
 \begin{pmatrix} 
 a_1& & & & & & -a_3&-a_5&\cdots & -a_{2p-1}\\
 & a_1 & & & & & -a_4&-a_6& \cdots & -a_{2p}\\
& & \ddots & & & &  a_2 & & &   \\
& & & a_1 &  & &  &a_2& &  \\
&  & & & a_1 & &  & &\ddots & \\
&  & & & & a_1 &  & &&a_2\\
-a_2& 0&-a_3&-a_5&\cdots & -a_{2p-1}& 0&\cdots & & 0\\
0& -a_2&-a_4&-a_6&\cdots & -a_{2p}& 0&\cdots & & 0
 \end{pmatrix}
\ene
This is a $2(\bbb-3)$-dimensional space of $\bbb\times 2\bbb-6$ matrices of bounded rank $\bbb-2$
generalizing  Case (III).

The example  Case (IV) generalizes to a $3q$-dimensional space of
$[2q+\binom q2]\times[2q+2]$ matrices of bounded rank $2q$.

More generally, going beyond corank two using the same blow-up, with $X$ a $k\times k$ matrix
of linear forms
$$
\begin{pmatrix}
a_1\Id_k & & -X\\
& a_1\Id_k & a_2\Id_k\\
a_2\Id_k & X& 0\end{pmatrix}
$$
gives a $k^2+2$ dimensional space of $3k\times 3k$ matrices of bounded rank $2k$.
 L. Manivel points out  that this suggests (after permuting the blocks):
$$
\begin{pmatrix}
a_1\Id_k &a_2\Id_k  &a_3\Id_k\\
X&  & &(a_2-a_3)\Id_k\\
  & X&  & (a_3-a_1)\Id_k\\
    &  & X& (a_1-a_2)\Id_k
  \end{pmatrix}
$$
and more generally $p$ blocks with $X$ and $(a_1\hd a_k)$,
 $(b_1\hd b_k)$ linear forms such that
$$(a_1\hd a_k)\begin{pmatrix} b_1\\ \vdots\\ b_k\end{pmatrix} =0.
$$

 There is a unique up to isomorphism concise tensor in $\BC^3\ot \BC^3\ot \BC^3$ that
 is both of minimal border rank and gives rise to a space of bounded rank \cite{MR3239293}, the space is
 $$
 \begin{pmatrix} a_1&a_2&a_3\\ & a_1& \\ &a_3&\end{pmatrix}
 $$
 (the tensor is $1_B$ and $1_C$-generic, and  $1_A$-degenerate).
 Applying the construction of Proposition \ref{blowupprop} with the other tensor $T_{skewcw,2}$, yields
 a bounded rank $6$ space in $\BC^9\ot \BC^9$.
 This space is compression, but it may be useful for Strassen's laser method. 
 It also shows that even if the space of $k\times k$ matrices is of bounded rank less than $k$, one
 can still obtain a bounded rank $kr$ space with the construction.

    \bibliographystyle{amsplain}

\bibliography{Lmatrix}

\end{document}